\documentclass{elsart}
\usepackage{epsfig}
\usepackage{amsfonts}
\usepackage{amssymb}

\def\Prob{{\hbox{\rm  Prob}}}

\def\inter{\hbox{\rm int}}
\def\bZd{\Bbb Z^d}
\def\bZ{\Bbb Z}
\def\bN{\Bbb N}
\def\bRd{\Bbb R^d}
\def\b1{\mathbf 1}
\def\cal{\mathbf}
\newcommand{\ese}{\end{eqnarray*}}
\newcommand{\bse}{\begin{eqnarray*}}
\begin{document}
\begin{frontmatter}
\title{Nonparametric denoising Signals of Unknown Local Structure, II:  Nonparametric Regression Estimation}
\author[Tolya]{Anatoli Juditsky\corauthref{cor}},
\corauth[cor]{Corresponding author.}
\ead{anatoli.juditsky@imag.fr}
\author[Arik]{Arkadi Nemirovski}
\ead{nemirovs@isye.gatech.edu}
\address[Tolya]{LJK, B.P. 53, 38041 Grenoble Cedex 9, France}
\address[Arik]{ISyE,
Georgia Institute of Technology,
765 Ferst Drive,
 Atlanta GA 30332-0205 USA}
\begin{abstract}
We consider the problem of recovering of continuous multi-dimensional functions
$f$
from the noisy observations over the regular grid $m^{-1}\bZd$, $m\in \bN_*$.
Our focus is at the adaptive estimation in the case
when the function can be well recovered using a
linear filter, which can depend on the unknown
function itself.
In the companion paper \cite{filters1} we have shown in the case when there exists an adapted time-invariant filter, which locally recovers ``well'' the unknown signal, there is a numerically efficient construction of an adaptive filter which recovers the signals ``almost as well''.
In the current paper we study the application of the proposed estimation techniques in the non-parametric regression setting. Namely, we propose an adaptive estimation procedure for ``locally well-filtered" signals (some typical examples being smooth signals,  modulated smooth signals and harmonic functions) and show that the rate of recovery of such signals in the $\ell_p$-norm on the grid is essentially the same as that rate for regular signals with nonhomogeneous smoothness.
\end{abstract}
\end{frontmatter}
\section{Introduction}
 Let ${\cal F}=(\Omega,\Sigma,P)$ be a probability space. We
consider the problem of recovering unknown complex-valued random field
$(s_\tau=s_\tau(\xi))_{{\tau\in\bZd\atop\xi\in\Omega}}$ over
$\bZd$ from noisy observations
\begin{equation}\label{eq:mod}
y_{\tau}=s_{\tau}+e_{\tau}.
\end{equation}
 We assume that
the field  $(e_{\tau})$ of observation noises is independent of
$(s_{\tau})$ and is of the form $e_{\tau}=\sigma \epsilon_{\tau}$,
where $(\epsilon_{\tau})$ are independent of each other {\sl
standard} Gaussian complex-valued variables; the adjective
``standard'' means that $\Re(\epsilon_{\tau})$,
$\Im(\epsilon_{\tau})$ are independent of each other ${\cal
N}(0,1)$ random variables.

We suppose that the observations
(\ref{eq:mod}) come from a function (``signal'') $f$ of continuous
argument (which we assume to vary in the $d$-dimensional unit cube
$[0,1]^d$); this function is observed in noise along an $n$-point
equidistant grid in $[0,1]^d$, and the problem is to recover $f$
via these observations. This problem fits the framework of nonparametric regression estimation with a ``traditional setting'' as follows:
\begin{description}
\item[]A. The objective is to recover an unknown  smooth function $f:\,[0,1]^d\to {\Bbb R}$,
which is sampled on the observation grid $\Gamma_n=\{x_\tau=m^{-1}\tau: \,0\le \tau_1,...,\tau_d\le m\}$
with $(m+1)^d=n$, so that $s_\tau=f(x_\tau)$.
The error of recovery is measured with some functional norm (or a semi-norm) $\|\cdot\|$ on $[0,1]^d$, and the risk of recovery $\widehat{f}$ of $f$ is the expectation
$E_f\|\widehat{f}-f\|^2$;\\
\item[]B. The estimation routines are aimed at
recovering {\sl smooth} signals, and their quality is measured by
their maximal risks,
 the maximum being taken over $f$ running through natural families
 of smooth signals, e.g., H\"older or Sobolev balls;\\
\item[]C. The focus is on the asymptotic, as the volume of observations
$n$ goes to infinity, behavior of the estimation routines, with
emphasis on  {\sl asymptotically minimax (nearly) optimal
estimates} -- those which attain (nearly) best possible rates of
convergence of the risks to 0 as the observation sample size $n\to\infty$.
\end{description}
Initially, the research was focused on recovering smooth signals
with {\sl a priori known smoothness parameters} and the estimation
routines were tuned to these parameters (see, e.g.,
\cite{bib:ibkh80,MS,bib:stone80,bib:ibkh81,Birge,bib:nem85,Wa,Ha,Ros,Go,KT}).
Later on, there was a significant research on {\sl adaptive
estimation}. Adaptive estimation methods are free of a priori assumptions on
the smoothness parameters of the signal to be recovered, and the
primary goal is to develop the routines which exhibit asymptotically
(nearly) optimal behavior on a wide variety of families of smooth functions
(cf. \cite{bib:pief84,bib:lep90,bib:lep91,bib:lepmamspok94,bib:donjohn92b,bib:donjohn92a,bib:donjohnker93,bib:jud94,BirgeMass,Don8,GN1}).
For a more compete overview of results on smooth nonparametric
regression estimation see, for instance, \cite{SF}.\footnotemark[3]
\footnotetext[3]{Our ``brief outline'' of adaptive approach to nonparametric
regression would be severely incomplete without mentioning a novel
approach aimed at recovering {\sl nonsmooth} signals possessing
{\sl sparse representations} in properly constructed functional
systems
\cite{Don11,Don12,Don4,Don5,Don6,Don9,Don10,Don7,Don22,Don2,Don1,Elad}.
This promising approach is completely beyond the scope of our
paper.}
\\
The traditional focus on recovering smooth signals ultimately
comes from the fact that such a signal {\sl locally} can be
well-approximated by a polynomial of a fixed order $r$, and such
a polynomial is an ``easy to estimate'' entity. Specifically, for
every integer $T\geq0$, the value of a polynomial $p$ at an
observation point $x_t$ can be recovered via $(2T+1)^d$
neighboring observations $\{x_\tau:|\tau_j-t_j|\leq T,1\leq j\leq
d\}$ ``at a parametric rate'' -- with the expected squared error
$C\sigma^2(2T+1)^{-d}$ which is inverse proportional to the amount
$(2T+1)^d$ of the observations used by the estimate. The
coefficient $C$ depends solely on the order $r$ and the dimensionality
$d$ of the polynomial. The corresponding estimate
$\widehat{p}(x_t)$ of $p(x_t)$ is pretty simple: it is given by
a ``time-invariant filter'', that is, by convolution of
observations with an appropriate discrete kernel
$q^{(T)}=(q^{(T)}_\tau)_{\tau\in\bZd}$ vanishing outside the box
${\cal O}_T=\{\tau\in\bZd:|\tau_j|\leq T,1\leq j\leq d\}$:
$$
\widehat{p}(x_t)=\sum\limits_{\tau\in{\cal O}_T} q^{(T)}_\tau
y_{t-\tau},
$$
then the estimation $\widehat{f}$ of $f(x_t)$ is taken as $\widehat{f}=\widehat{p}(x_t)$.
\par
Note that the kernel $q^{(T)}$ is readily given
by the degree $r$ of the approximating polynomial,  $T$ and dimension $d$.
The ``classical" adaptation routines takes care of choosing ``good" values of the approximation parameters
(namely, $T$ and $r$). On the other hand, the polynomial approximation ``mechanism"
is supposed to be fixed  once for ever. Thus, in those procedures the ``form" of the kernel is considered
as given in advance.
\par
In the companion paper \cite{filters1} (referred hereafter as Part I) we have introduced the notion of a {\em well-filtered signal}. In brief,
the signal $(s_\tau)_{\tau\in\bZd}$ is $T$-well-filtered for some $T\in \bN_+$ if there is a filter (kernel) $q=q^{(T)}\in \cal O_T$
which recovers $(s_\tau)$ in the box
$\{u:|u-t|\leq 3T\}$ with the mean square error comparable with
$\sigma T^{-d/2}$:%
%
\[
\max\limits_{u:|u-t|\leq 3T}E\left\{|s_u-\sum\limits_{\tau\in{\cal O}_T} q^{(T)}_\tau
y_{u-\tau}|^2\right\}\leq
O(\sigma^2T^{-d}).
\]
 The universe of
these signals is much wider than the one of smooth signals.
As we have seen in Part I that it contains, in particular, ``modulated smooth signals''
-- sums of a fixed number of products of  smooth functions and
multivariate harmonic oscillations of unknown (and arbitrarily
high) frequencies. We have shown in Part I that
whenever a discrete time signal (that is, a signal defined on
a regular discrete grid) is {\sl well-filtered}, we can recover this
signal at a ``nearly parametric'' rate {\sl without a priori
knowledge of the associated filter}. In other words, a well-filtered signal can be recovered {\sl on the observation grid}
 basically as well as if it were an algebraic polynomial
 of a given order.

We are about to demonstrate that the results of Part I on recovering well-filtered signals of unknown
structure can be applied to recovering nonparametric signals which
admit well-filtered {\sl local} approximations. Such an
extension has an unavoidable price -- now we cannot hope to
recover the signal well outside of the observation grid (a highly oscillating
 signal
 can merely
vanish on the observation grid and be arbitrarily large outside it). As a
result, {\sl in what follows we are interested in recovering the
signals along the observation grid only} and, consequently,
replace the error measures based on the functional norms on $[0,1]^d$ by their grid analogies. 
\par
  The estimates to be developed will be ``double adaptive'', that
is, adaptive with respect to both the unknown in advance
structures of well-filtered approximations of our signals and to
the unknown in advance ``approximation rate'' -- the dependence
between the size of a neighborhood of a point where the signal in
question is approximated and the quality of approximation in this
neighborhood. Note that in the case of smooth signals, this approximation
rate is exactly what is given by the smoothness parameters. The
results to follow can be seen as extensions of the results of
\cite{bib:nem92,GN2} (see also \cite{SF}) dealing with the
particular case of univariate signals satisfying differential
inequalities with unknown differential operators.


\section{Nonparametric regression problem}
We start with the formal description of the components of the nonparametric regression problem.
\par
\par

Let for $\tau\in \bZd$, $|\tau|=\max \{|\tau_1|,...,|\tau_d|\}$, and let $\tau\le m$ for some $a\in \bN$ denote $\tau_i\le m,\;i=1,...,d$.
Let $m$ be a positive integer, $n=(m+1)^d$, and let
$\Gamma_n=\{x=m^{-1}\alpha:\alpha\in \bZd,0\leq\alpha,|\alpha|\leq
m\}$.
 \par
 Let $C([0,1]^d)$ be the linear space of complex-valued fields over $[0,1]^d$.
 We associate with a signal $f\in C([0,1]^d)$ its
observations along $\Gamma_n$:
\begin{equation}\label{3.eq1}
y\equiv y^n_f(\epsilon)=\{y_\tau\equiv y_\tau^n(f,\epsilon) =
f(m^{-1}\tau)+e_\tau, e_\tau=\sigma
\epsilon_\tau\}_{0\le \tau\le m},
\end{equation}
where $\{\epsilon_\tau\}_{\tau\in\bZd}$ are independent standard
Gaussian complex-valued random noises. Our goal is to recover
$f\big|_{\Gamma_n}$ from observations (\ref{3.eq1}).
  In what
follows, we write
$$
f_\tau = f(m^{-1}\tau),\;\;\;\;\;{[\tau\in\bZd,m^{-1}\tau\in[0,1]^d]}
$$

 Below we use the following notations.
  For a set $B\subset[0,1]^d$, we denote by $\bZ(B)$ the set of all
$t\in\bZd$ such that $m^{-1}t\in B$.  We denote $\|\cdot\|_{q,B}$ the standard $L_p$-norm on $B$:
$$
\|g\|_{p,B}= \left(\,\int_{x\in B}|g(x)|^p dx\right)^{1/p},
$$
and $|g|_{q,B}$ its discrete analogy, so that
$$
|g|_{q,B} = m^{-d/q}\left(\sum\limits_{\tau\in
\bZ(B)}|g_\tau|^q\right)^{1/q}.
$$
We set
 $$\Gamma_n^o=\Gamma_n\cap(0,1)^n=\{m^{-1}t: t\in \bZd,
 t>0,|t|<m\}.$$
Let
$x=m^{-1}t\in\Gamma_n^o$. We say that a nonempty {\sl open} cube
$$
{B}_h(x) =\{u\mid\, |u_i-x_i|< h/2,\,\,i=1,...,d\}
$$
centered at $x$ is {\sl admissible} for $x$, if
$B_h(x)\subset[0,1]^n$. For such a cube, $T_h(x)$ denotes the
largest nonnegative integer $T$ such that $$\bZ(B)\supset
\{\tau\in\bZd: |\tau-t|\leq 4T\}.$$    For a cube $$B=\{x\in\bRd:
|x_i-c_i|\leq h/2,\,i=1,...,d\},$$  $D(B)=h$ stands for the edge of
$B$. For $\gamma\in(0,1)$ 
we denote   $$B_\gamma=\{x\in\bRd: |x_i-c_i|\leq\gamma
h/2,i=1,...,d\}$$ the $\gamma$-shrinkage of $B$ to the center
of $B$.
\subsection{Classes of locally well-filtered signals}
%
Recall that we say that a function on $[0,1]^d$ is smooth if it can be locally well-approximated by a polynomial. {\em Informally}, the
the definition below sais that  a continuous signal $f\in C([0,1])^d$ is {\em locally well-filtered}  if $f$ admits a good  {\em local} approximation by  a well-filtered discrete signal $\phi_\tau$ on $\Gamma_n$ (see Definition 1 of Section 2.1, Part I).
\begin{defn}\label{def:wfs}
Let
$B\subset[0,1]^d$ be a cube, $k$ be a positive integer,
$\rho\geq1$, $R\geq0$ be reals, and let $p\in(d,\infty]$.
The
collection $B$, $k$, $\rho$, $R$, $p$ specifies the family ${\cal
F}^{k,\rho,p}(B,R)$ of {\sl locally well-filtered on $B$}  signals
$f$ defined by the following requirements:
\\
(1) $f\in C([0,1]^d)$; \\
(2) There exists a nonnegative function $F\in L_p(B),\;\|F\|_{p,B}\leq R,$
such that for every
$x=m^{-1}t\in\Gamma_n\cap\inter B$ and for every admissible for
$x$ cube $B_h(x)$ contained in $B$ there exists a field $\phi\in
C(\bZd)$ such that $\phi\in {\cal S}^t_{3T_h(x)}(0,\rho,T_h(x))$ (where the set ${\cal S}^t_L(\theta,\rho,T)$ of $T$-well filtered signals is defined in Definition 1 of Part I)
and
\begin{equation}\label{dapp:main}
\forall\tau\in \bZ(B_h(x)): |\phi_\tau-f_\tau|\leq
h^{k-d/p}\|F\|_{p,B_h(x)}.
\end{equation}
\end{defn}
 In the sequel, we use for ${\cal F}^{k,\rho,p}(B;R)$
also the shortened notation ${\cal F}[\psi]$, where $\psi$ stands
for the collection of ``parameters'' $(k,\rho,p,B,R)$.
%
\par {\bf Remark} The motivating example of locally well-filtered signals is that of modulated smooth signals as follows. Let a cube
$B\subset[0,1]^d$, $p\in(d,\infty]$, positive integers $k,\nu$ and
a real $R\geq0$ be given. Consider a collection of $\nu$ functions
$g_1,...,g_\nu\in C([0,1]^d)$ which are $k$ times continuously
differentiable and satisfy the constraint
$$
\sum\limits_{\ell=1}^\nu\|D^k g_\ell\|_{p,B}\leq R. $$
Let
$\omega(\ell)\in \bRd$, and let $$ f(x)=\sum\limits_{\ell=1}^\nu
g_\ell(x)\exp\{i\omega^T(\ell)x\}. $$
By the standard argument
\cite{myref1}, whenever $x=m^{-1}t\in\Gamma_n\cap\inter B$ and
$B_h(x)$ is admissible for $x$, the Taylor polynomial
$\Phi_\ell^x(\cdot)$, of order $k-1$, taken at $x$, of $f_\ell$
satisfies the inequality
$$ u\in B_h(x)\Rightarrow
|\Phi_\ell^x(u)-f_\ell(u)|\leq
c_1h^{k-d/p}\|F_\ell\|_{p,B_h(x)},\;\;\mbox{where} \;\; F_\ell(u)=|D^kf_\ell(u)|
$$
(here and in what follows, $c_i$ are positive constants depending
solely on $d$, $k$ and $\nu$). It follows that if $
\Phi(u)=\sum\limits_{\ell=1}^\nu
\Phi_\ell^{x}(u)\exp\{i\omega^T(\ell)u\} $ then
\begin{equation}\label{dapp:close}
\begin{array}{l}
u\in B_h(x)\Rightarrow |\Phi(u)-f(u)|\leq
h^{k-d/p}\|F\|_{p,B_h(x)},\\
 F=c_2\sum\limits_{\ell=1}^\nu F_\ell
\qquad [\Rightarrow \|F\|_{p,B}\leq c_3R].\\
\end{array}
\end{equation}
 Now observe that the exponential polynomial
$ \phi(\tau)=\Phi(m^{-1}\tau) $ belongs to ${\cal
S}^t_L(0,c_4,T)$ for any $0\le T\le L\le \infty$ (Proposition 10 of Part I). Combining this
fact with (\ref{dapp:close}), we conclude that $f\in {\cal
F}^{k,\rho(\nu,k,d),p}(B,c(\nu,k,d)R).$
\subsection{Accuracy measures}
\label{accmeas} Let us fix $\gamma\in(0,1)$ and $q\in [1,\infty]$. Given an estimate
$\widehat{f}_n$ of the restriction $f\big|_{\Gamma_n}$ of $f$ on the grid ${\Gamma_n}$, based on observations
(\ref{3.eq1}) (i.e., a Borel real-valued function of $x\in
\Gamma_n$ and $y\in{\Bbb C}^n$) and $\psi=(k,\rho,p,B,R)$, let us
characterize the quality of the estimate on the set ${\cal
F}[\psi]$ by the worst-case risks
$$
\widehat{{\cal R}}_q\left(\widehat{f}_n;{\cal F}[\psi]\right)=
\sup_{f\in {\cal F}[\psi]} \left(E
\left\{\left|\widehat{f}_n(\cdot;y_f(\epsilon))-f\big|_{\Gamma_n}(\cdot)\right|_{
q,B_\gamma}^2\right\}\right)^{1/2}.
$$

\section{Estimator construction} \label{3.theestimate} The recovering
routine we are about to build is aimed at estimating functions
from classes ${\cal F}^{k,\rho,p}(B,R)$ with {\sl unknown in
advance parameters $k,\rho,p,B,R$}. The only design parameters of
the routine is an a priori upper bound $\mu$ on the parameter
$\rho$ and a $\gamma\in(0,1)$.
\subsection{Preliminaries} From now on, we denote by
$\Theta=\Theta_{(n)}$ the deterministic function of observation
noises defined as follows. For every cube $B\subset[0,1]^d$ with
vertices in $\Gamma_n$, we consider the discrete Fourier transform of the
observation noises reduced to $B\cap\Gamma_n$, and take the
maximum of modules of the resulting Fourier coefficients, let it
be denoted $\theta_B(e)$. By definition,
$$
\Theta\equiv \Theta_{(n)}=\sigma^{-1}\max\limits_B \theta_B(e),
$$
where the maximum is taken over all cubes $B$ of the indicated
type. By the origin of $\Theta_{(n)}$, due to the classical results on maxima of Gaussian processes (cf also Lemma 15 of Part I), we have
\begin{equation}\label{3.beq18}
\forall w\ge1:\quad \Prob\left\{\Theta_{(n)}> w\sqrt{\ln n}
\right\} \le \exp\left\{-{c{w^2\ln n}\over 2}\right\},
\end{equation}
where $c>0$ depends solely on $d$.
\subsection{Building blocks: window estimates} To recover a
signal $f$ via $n=m^d$ observations (\ref{3.eq1}), we use
point-wise window estimates of $f$ defined as follows.
\par
Let us fix
a point $x=m^{-1}t\in\Gamma_n^o$; our goal is to build an estimate
of $f(x)$. Let $B_h(x)$ be an admissible window for $x$. We
associate with this window an estimate
$\widehat{f}^h_n=\widehat{f}^h_n(x;y^n_f(\epsilon))$ of $f(x)$
defined as follows. If the window is ``very small'', specifically,
$h\leq m^{-1}$, so that $x$ is the only point from the observation
grid $\Gamma_n$ in $B_h(x)$, we set $T_h(x)=0$ and
$\widehat{f}^h_n=y_t$. For a larger window, we choose the largest
nonnegative integer $T=T_h(x)$ such that
$$\bZ(B_h(x))\supset\{\tau:|\tau-t|\leq 4T\}$$ and apply Algorithm
A of Part I to build the estimate of $f_t=f(x)$, the design parameters of the algorithm being $(\mu ,T_h(x))$. Let the resulting estimate be
denoted by $\widehat{f}^h_n=\widehat{f}^h_n(x;y^n_f(\epsilon))$.
\\  To characterize the quality of the estimate
$\widehat{f}^h_n=\widehat{f}^h_n(x;y^n_f(\epsilon))$, let us set
$$
\Phi_\mu(f,{B}_h(x))=\min\limits_p\left\{\max\limits_{\tau\in
\bZ(B_h(x))}|p_\tau-f_\tau|: p \in {\cal
S}^t_{3T_h(x)}(0,\mu,T_h(x))\right\}.
$$
\begin{lem}\label{dapp:lemma1} One has
\begin{equation}\label{dapp:em10}
(f_\tau)\in {\cal S}^t_{3T_h(x)}(\theta,\mu,T_h(x)),\quad \theta=
{\Phi_\mu(f,B_h(x))(1+\mu)\over(2T+1)^{d/2}}.
\end{equation}
\end{lem}
Assuming that $h>m^{-1}$ and combining (\ref{dapp:em10}) with
the result of Theorem 4 of Part I we come to the following upper bound on the
error of estimating $f(x)$ by the estimate
$\widehat{f}^h_n(x;\cdot)$:
\begin{equation}
\label{3.eq5} |f(x)-\widehat{f}^h_n(x;y_f(\epsilon))|\leq
C_1\left[\Phi_\mu(f,B_h(x))+
{\sigma\over{\sqrt{nh^d}}}\Theta_{(n)}\right]
\end{equation}
(note that $(2T_h(x)+1)^{-d/2}\leq C_0(nh^d)^{-1/2}$). 
For evident reasons  (\ref{3.eq5}) holds true for
``very small windows'' (those with $h\leq m^{-1}$) as well.
\subsection{The adaptive estimate}
We are about to ``aggregate'' the window estimates
$\widehat{f}^h_n$ into an {\sl adaptive} estimate, applying
Lepskii's adaptation scheme in the
same fashion as in \cite{bib:lepmamspok94,GN1,GN2}.
\\
 Let us fix a ``safety
factor'' $\omega$ in such a way that the event
$\Theta_{(n)}>\omega\sqrt{\ln n}$ is ``highly un-probable'',
namely,
\begin{equation}
\label{3.beq18a} \Prob\left\{\Theta_{(n)}>\omega\sqrt{\ln
n}\right\}\le n^{-4(\mu+1)};
\end{equation}
by (\ref{3.beq18}), the required $\omega$ may be chosen as a
function of $\mu,d$ only. We are to describe the basic blocks of the construction of the adaptive estimate.\\
{\bf ``Good'' realizations of noise.} Let us define the set of
``good realizations of noise'' as
\begin{equation}
\label{3.eq9} \Xi_n=\{\epsilon\mid\, \Theta_{(n)}\le
\omega\sqrt{\ln n}\}.
\end{equation}
Now (\ref{3.eq5}) implies the ``conditional'' error bound
\begin{equation}
\label{3.eq7}\begin{array}{l} \epsilon\in\Xi_n\Rightarrow
|f(x)-\widehat{f}^h_n(x;y_f(\epsilon))|\leq
C_1\left[\Phi_\mu(f,{B}_h(x))+
S_n(h)\right],\\[3mm]
S_n(h)=\displaystyle{{\sigma\over{\sqrt{nh^d}}}\omega\sqrt{\ln n}}.\\
\end{array}
\end{equation}
Observe that {\sl as $h$ grows, the
``deterministic term'' $\Phi_\mu(f,{B}_h(x))$ does not decrease,
while the ``stochastic term'' $S_n(h)$ decreases.}
\par\noindent
{\bf The ``ideal'' window.} Let us define the {\sl ideal} window
${B}_*(x)$ as the largest admissible window for which the
stochastic term dominates the deterministic one:
\begin{equation}
\label{3.eq8}
\begin{array}
{rcl}
{B}_*(x)&=&{B}_{h_*(x)}(x),\\
h_*(x)&=&\max\{h\mid\,h>0, {B}_h(x)\subset[0,1]^d,
\Phi_\mu(f,{B}_h(x))\le S_n(h)\}.\\
\end{array}
\end{equation}
Note that such a window does exist, since $S_n(h)\to\infty$ as
$h\to+0$.   Besides this, since the cubes $B_h(x)$ are open, the
quantity $\Phi_\mu(f,{B}_h(x))$ is continuous from the left, so
that
\begin{equation}\label{miutar}
0<h\leq h_*(x)\Rightarrow \Phi_\mu(f,{B}_h(x))\le S_n(h).
\end{equation}
Thus, the ideal window ${B}_*(x)$ is well-defined for every $x$
possessing admissible windows, i.e., for every
$x=\Gamma_n^o=\{m^{-1}t: t\in\bZd, 0<t, |t|<m\}$. \\
{\bf Normal windows.} Assume that $\epsilon\in \Xi_n$.  Then the
errors of all estimates $\widehat{f}^h_n(x;y)$ associated with
admissible windows smaller than the ideal one are dominated by the
corresponding stochastic terms:
\begin{equation}
\label{3.eq10} \epsilon\in\Xi_n,0<h\leq h_*(x) \Rightarrow
|f(x)-\widehat{f}^h_n(x;y_f(\epsilon))| \le 2C_1S_n(h)
\end{equation}
(by (\ref{3.eq7}) and (\ref{miutar})). Let us fix an
$\epsilon\in\Xi_n$ (and thus -- a realization $y$ of the
observations) and let us call an admissible for $x$ window
${B}_h(x)$ {\sl normal}, if the associated estimate
$\widehat{f}^h_n(x;y)$ differs from every estimate associated with
a smaller  window by no more than $4C_1$ times the stochastic term
of the latter estimate, i.e.
\begin{equation}
\label{3.eq12}
\begin{array}
{c}
\hbox{Window ${B}_h(x)$ is normal}\\
\Updownarrow\\
\left\{\begin{array} {l}
\hbox{${B}_h(x)$ is admissible}\\
\forall h', 0<h'\leq h:\quad
|\widehat{f}^{h'}_n(x;y)-\widehat{f}^h_n(x;y)|\le 4C_1S_n(h')\quad
[y=y_f(\epsilon)]\\
\end{array}\right.\\
\end{array}
\end{equation}
Note that if $x\in\Gamma_n^o$, then $x$ possesses a normal window,
specifically, the window $B_{m^{-1}}(x)$. Indeed, this window
contains a single observation point, namely, $x$ itself, so that
the corresponding estimate, same as every estimate corresponding
to a smaller window, by construction coincides with the
observation at $x$, so that all the estimates
$\widehat{f}^{h'}_n(x;y)$, $0<h'\leq m^{-1}$, are the same. Note
also that (\ref{3.eq10}) implies that
\begin{quote}
(!) {\sl If $\epsilon\in\Xi_n$, then the ideal window ${B}_*(x)$
is normal}.
\end{quote}
{\bf The adaptive estimate $\widehat{f}_n(x;y)$.} The property of
an admissible window to be normal is ``observable'' -- given
observations $y$, we can say whether a given window is or is not
normal. Besides this, it is clear that among all normal windows
there exists the largest one ${B}^+(x)={B}_{h^+(x)}(x)$. {\sl The
adaptive estimate $\widehat{f}_n(x;y)$ is exactly the window
estimate associated with the window $B^+(x)$.}   Note that from
(!) it follows that
\begin{quote}
(!!) {\sl If $\epsilon\in \Xi_n$, then the largest normal window
${B}^+(x)$ contains the ideal window ${B}_*(x)$.}
\end{quote}
By definition of a normal window, under the premise of (!!) we
have
$$
|\widehat{f}^{h^+(x)}_n(x;y)-\widehat{f}^{h_*(x)}_n(x;y)|\le
4C_1S_n(h_*(x)),
$$
and we come to the conclusion as follows:
\\
(*) {\sl If $\epsilon\in \Xi_n$, then the error of the estimate $
\widehat{f}_n(x;y)\equiv\widehat{f}^{h^+(x)}_n(x;y) $ is dominated
by the error bound {\rm (\ref{3.eq7})} associated with the ideal
window:}
\begin{equation}
\label{3.eq15} \epsilon\in\Xi_n\Rightarrow
|\widehat{f}_n(x;y)-f(x)|\le5C_1\left[\Phi_\mu(f,{B}_{h_*(x)}(x))+S_n(h_*(x))\right].
\end{equation}
Thus, the estimate $\widehat{f}_n(\cdot;\cdot)$ -- which is based
solely on observations and does not require any a priori knowledge
of the ``parameters of well-filterability of $f$'' -- possesses
basically the same accuracy as the ``ideal'' estimate associated
with the ideal window (provided, of course, that the realization
of noises is not ``pathological": $\epsilon\in\Xi_n$).
\\
Note that the adaptive estimate $\widehat{f}_n(x;y)$ we have built
 depends solely on ``design parameters'' $\mu$, $\gamma$ (recall that $C_1$ depends on $\mu,\gamma$),
 the volume of
 observations $n$ and the dimension $d$.
\section{Main result} Our main result is as follows:
\begin{thm}\label{3.themain} Let $\gamma\in(0,1)$, $\mu\ge1$ be an integer,
let ${\cal F}={\cal F}^{k,\rho,p}(B;R)$ be a family of locally
well-filtered signals associated with a cube $B\subset[0,1]^d$
with $mD(B)\geq1$, $\rho\leq\mu$ and $p>d$. For properly chosen
$P\ge1$ depending solely on $\mu,d,p,\gamma$ and nonincreasing in
$p>d$ the following statement holds true:
\par
Suppose that the volume $n=m^d$ of observations {\rm (\ref{3.eq1})} is large
enough, namely,
\begin{equation}
\label{3.large}  P^{-1}n^{{{2kp+d(p-2)}\over{2dp}}}\ge
{R\over{\sigma}}\sqrt{ n\over {\ln n}} \ge P
[D(B)]^{-{{2kp+d(p-2)}\over 2p}}
\end{equation}
where
$D(B)$ is the edge of the cube $B$.
\par Then for every
$q\in[1,\infty]$ the worst case, with respect to ${\cal F}$,
$q$-risk of the adaptive estimate $\widehat{f}_n(\cdot,\cdot)$
associated with the parameter $\mu$ can be bounded as follows:
\begin{eqnarray}
\label{3.main}
\widehat{{\cal R}}_q\left(\widehat{f}_n;{\cal F
}\right)&=& \sup\limits_{f\in {\cal F}} \left(E \left\{\left|\widehat{f}_n(\cdot;y_f(\epsilon))-f(\cdot)\right|_{q,B_\gamma}^2\right\}\right)^{1/2}\\
\nonumber
&\le&PR\left({\sigma^2 \ln n\over
R^2n}\right)^{\beta(p,k,d,q)}
[D(B)]^{d\lambda(p,k,d,q)},
\end{eqnarray}
where
\begin{eqnarray*}
\beta(p,k,d,q)&=&\left\{\begin{array}{lcl}
{k\over{2k+d}},&\mbox{when}&{q}\le
{(2k+d)p\over d},\\
 {{k+d\left({1\over q}-{1\over p}\right)}\over{2k+d-{2d\over p}}},&\mbox{when}&{q}>
{(2k+d)p\over d},
\end{array}\right.\\
\lambda(p,k,d,q)&=&\left\{\begin{array}{lcl}
{1\over q}-{d\over{(2k+d)p}},&\mbox{when}&{q}\le
{(2k+d)p\over d},\\
0,&\mbox{when}&{q}>
{(2k+d)p\over d}
\end{array}\right.\end{eqnarray*}
(recall that here $B_\gamma$ is the concentric to $B$ $\gamma$ times smaller cube).
\end{thm}
Note that the rates of convergence to 0, as $n\to \infty$, of the
risks $\widehat{{\cal R}}_q\left(\widehat{f}_n;{\cal F }\right)$
of our adaptive estimate  on the families ${\cal F}={\cal
F}^{k,\rho,p}(B;R)$ are exactly the same as those stated by
Theorem 3 from \cite{bib:nem85} (see also \cite{bib:lepmamspok94,bib:donjohnker93,GN1,SF}) in the case of recovering
non-parametric {\sl smooth} regression functions from
Sobolev balls. It is well-known  that in the
smooth case the latter rates are  optimal in order, up to
logarithmic in $n$ factors. Since the families of locally
well-filtered signals are much wider than local Sobolev balls
(smooth signals are trivial examples of modulated smooth
signals!), it follows that the rates of convergence stated by
Theorem \ref{3.themain} also are nearly optimal.
\section{Simulation examples}
\label{sectnum}
In this section we present the results of a small simulation study of the adaptive filtering algorithm applied to the 2-dimensional de-noising problem. The simulation setting is as follows:
we consider real-valued signals
\[
y_\tau=s_\tau+e_{\tau},\;\;\tau=(\tau_1,\,\tau_2)\in \{1,...,m\}^2,
\]
$e_{(1,1)},...,\,e_{(m,m)}$ being  independent standard Gaussian random variables.
The problem is to estimate, given observations $(y_\tau)$, the values of the signal $(f_{x_\tau})$ on the grid
$\Gamma_m=\{m^{-1}\tau,\;1\le \tau_1, \tau_2\le m\}$, and $f({x_\tau})=s_\tau$. The value $m=128$ is common to all experiments.

We consider signals which are sums of three harmonic components:
\[
s_\tau=\alpha[\sin(m^{-1}\omega_1^T\tau+\theta_1)+\sin(m^{-1}\omega_2^T\tau+\theta_2)+\sin(m^{-1}\omega_3^T\tau+\theta_3)];
\]
the frequencies $\omega_i$ and the phase shifts $\theta_i$, $i=1,...,3$ are drawn randomly from the uniform distribution over, respectively, $[0,\omega_{\max}]^3$ and $[0,1]^3$ and the coefficient $\alpha$ is chosen to obtain the signal-to-noise ratio equal to one.
\par
\begin{figure}
$$\begin{array}{cc}
\epsfxsize=190pt\epsfysize=190pt\epsffile{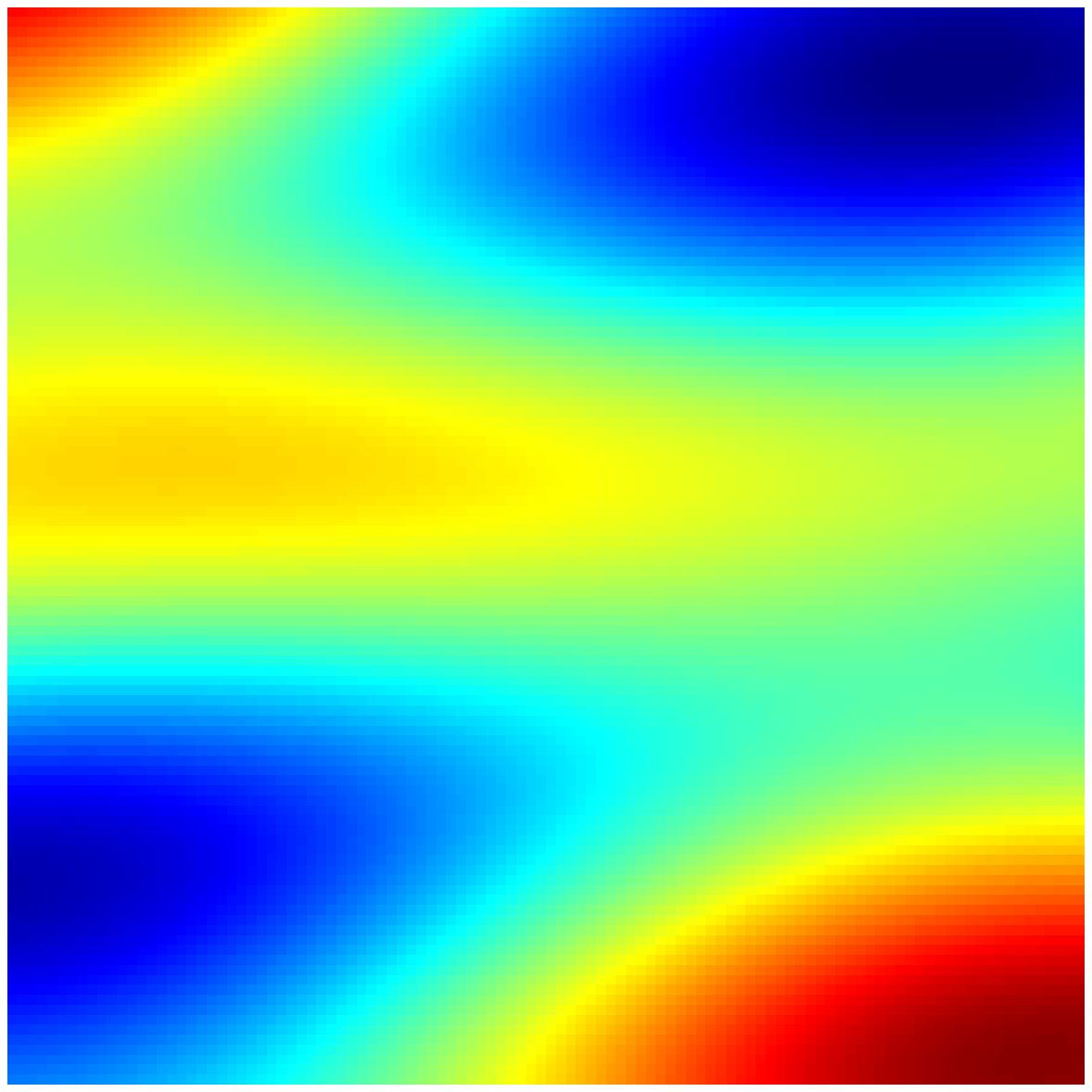}&
\epsfxsize=190pt\epsfysize=190pt\epsffile{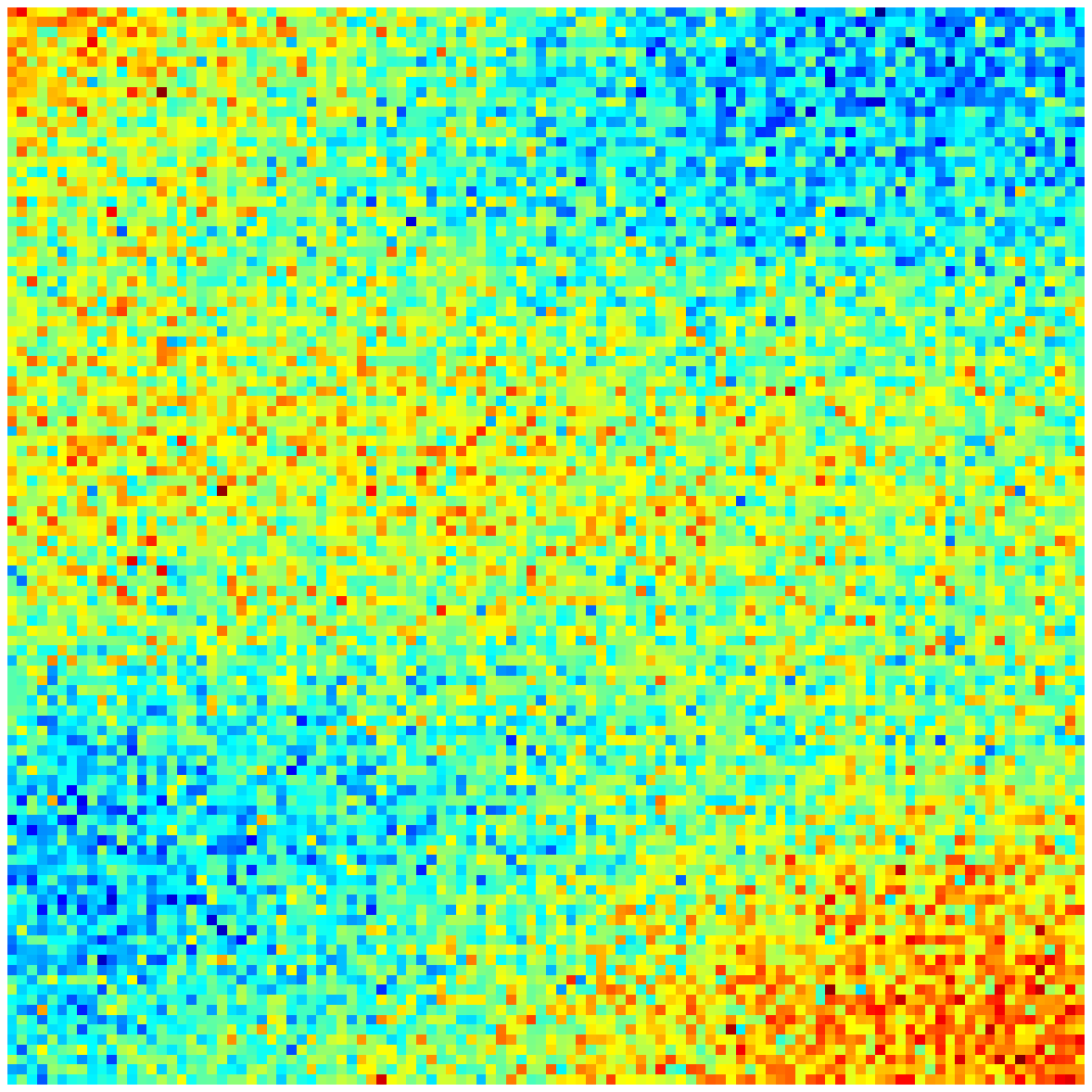}\\
\hbox{True Image}&\hbox{Observation} \\
\epsfxsize=190pt\epsfysize=190pt\epsffile{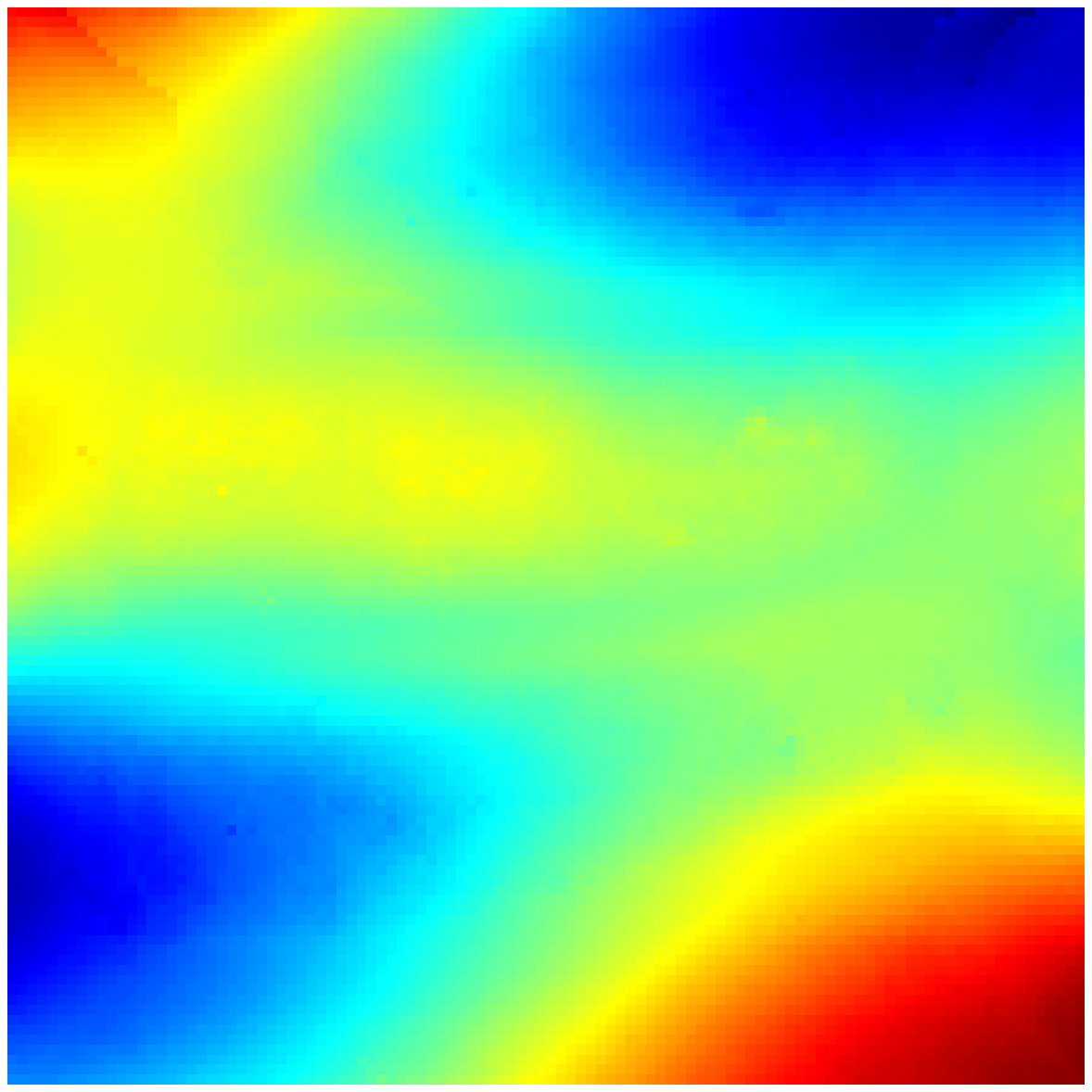}&
\epsfxsize=190pt\epsfysize=190pt\epsffile{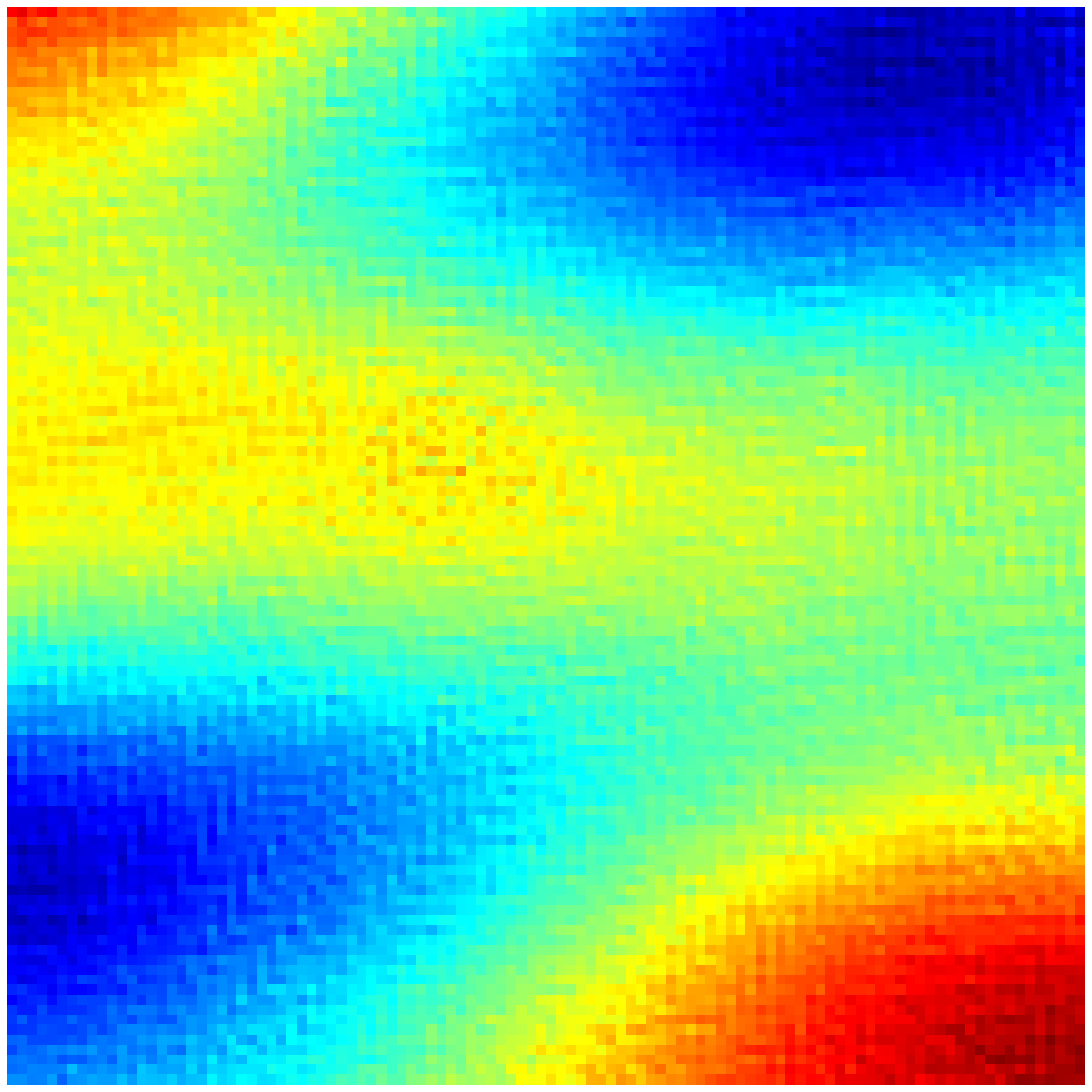}\\
\hbox{Standard recovery}&\hbox{Adaptive recovery}
\end{array} $$
\caption{Recovery for $\omega_{\max}=2.0$}\label{pic:1}\end{figure}
\begin{figure}
$$\begin{array}{cc}
\epsfxsize=190pt\epsfysize=190pt\epsffile{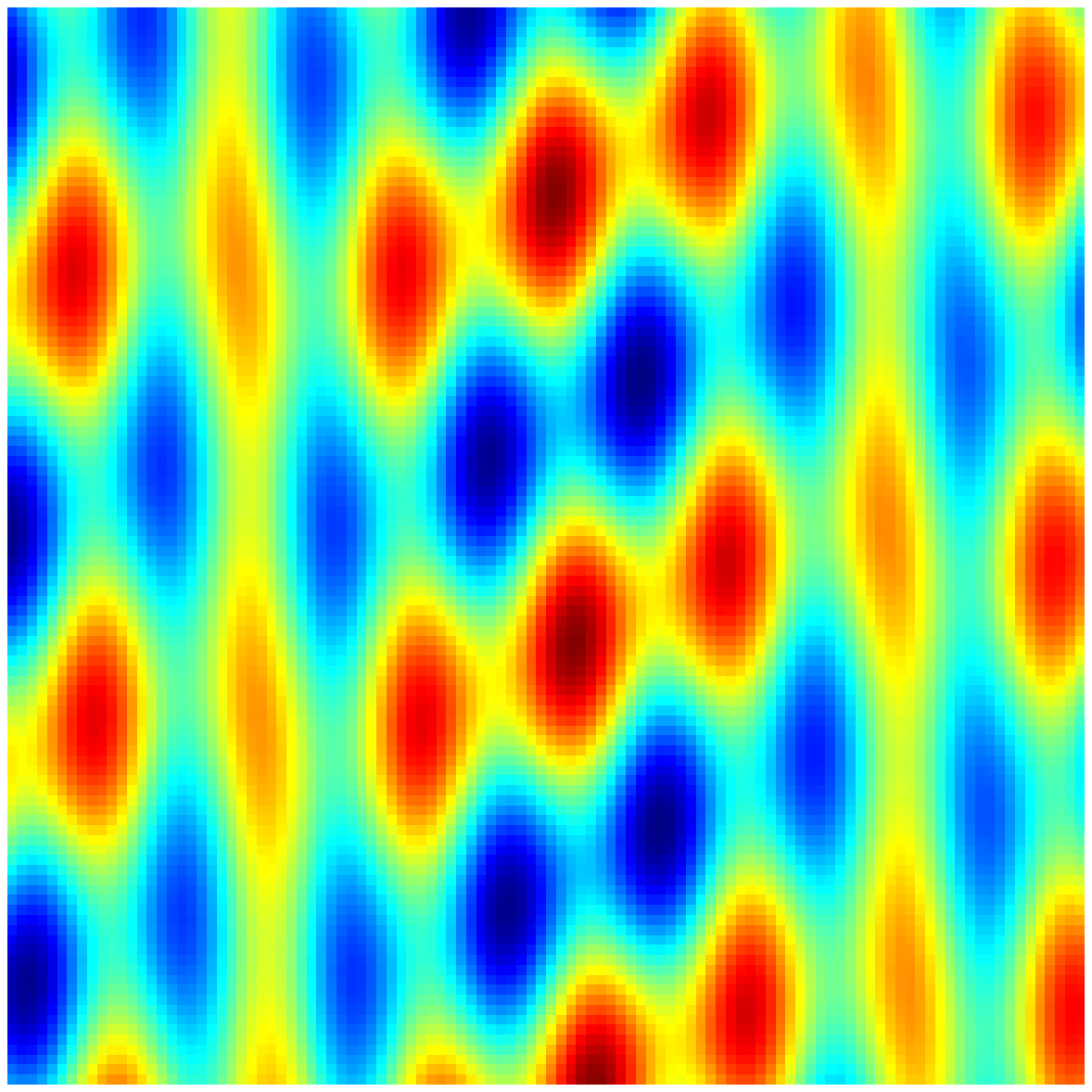}&
\epsfxsize=190pt\epsfysize=190pt\epsffile{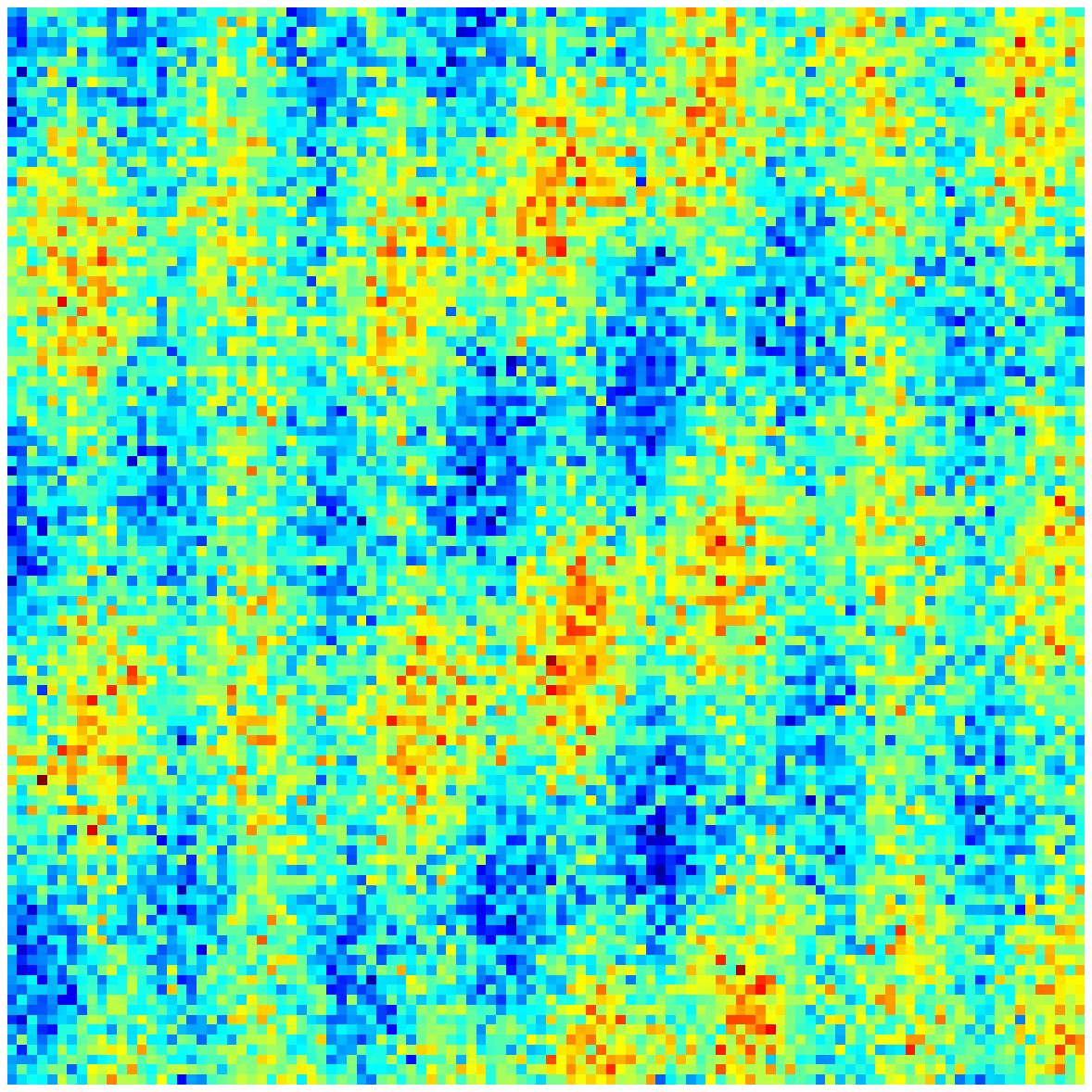}\\
\hbox{True Image}&\hbox{Observation} \\
\epsfxsize=190pt\epsfysize=190pt\epsffile{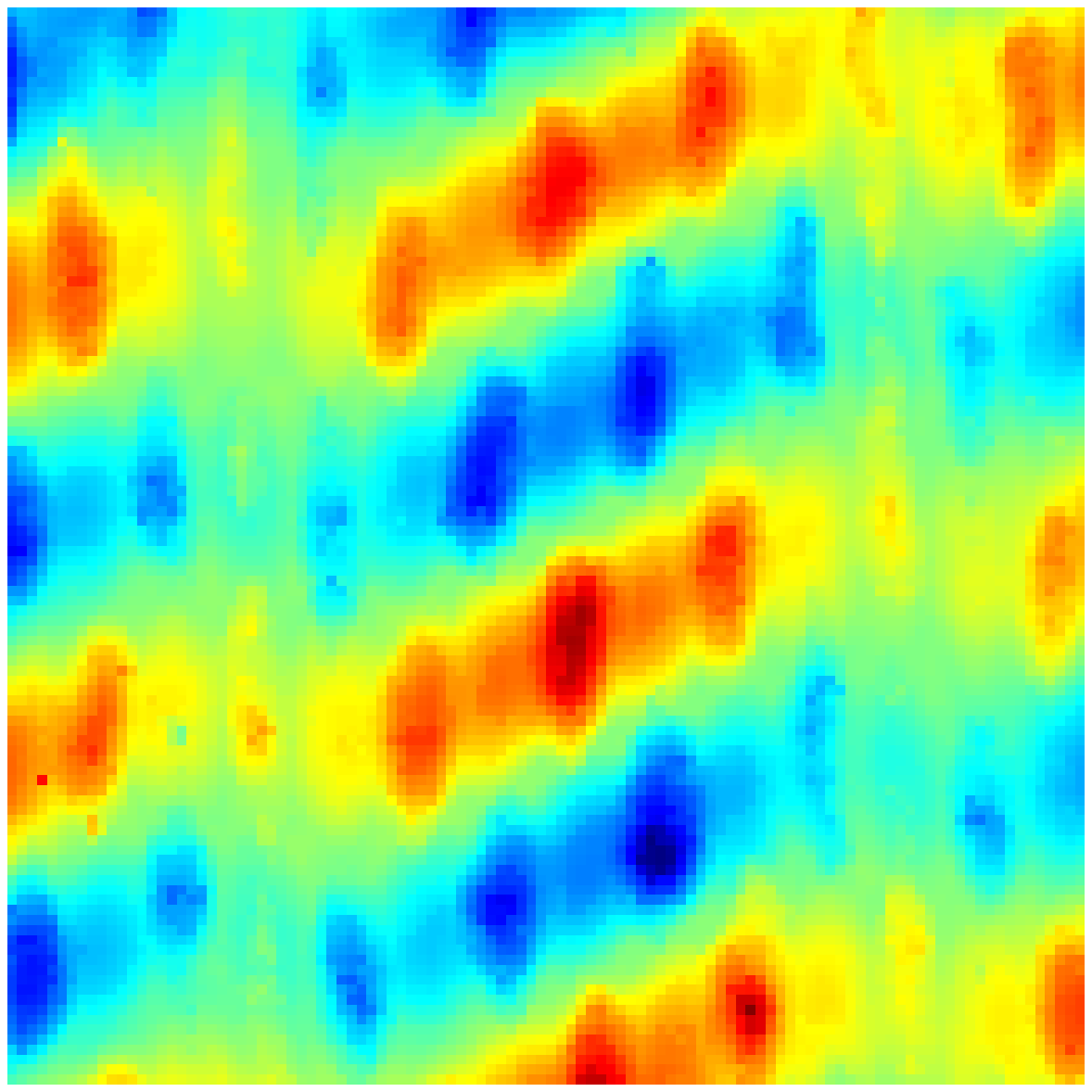}&
\epsfxsize=190pt\epsfysize=190pt\epsffile{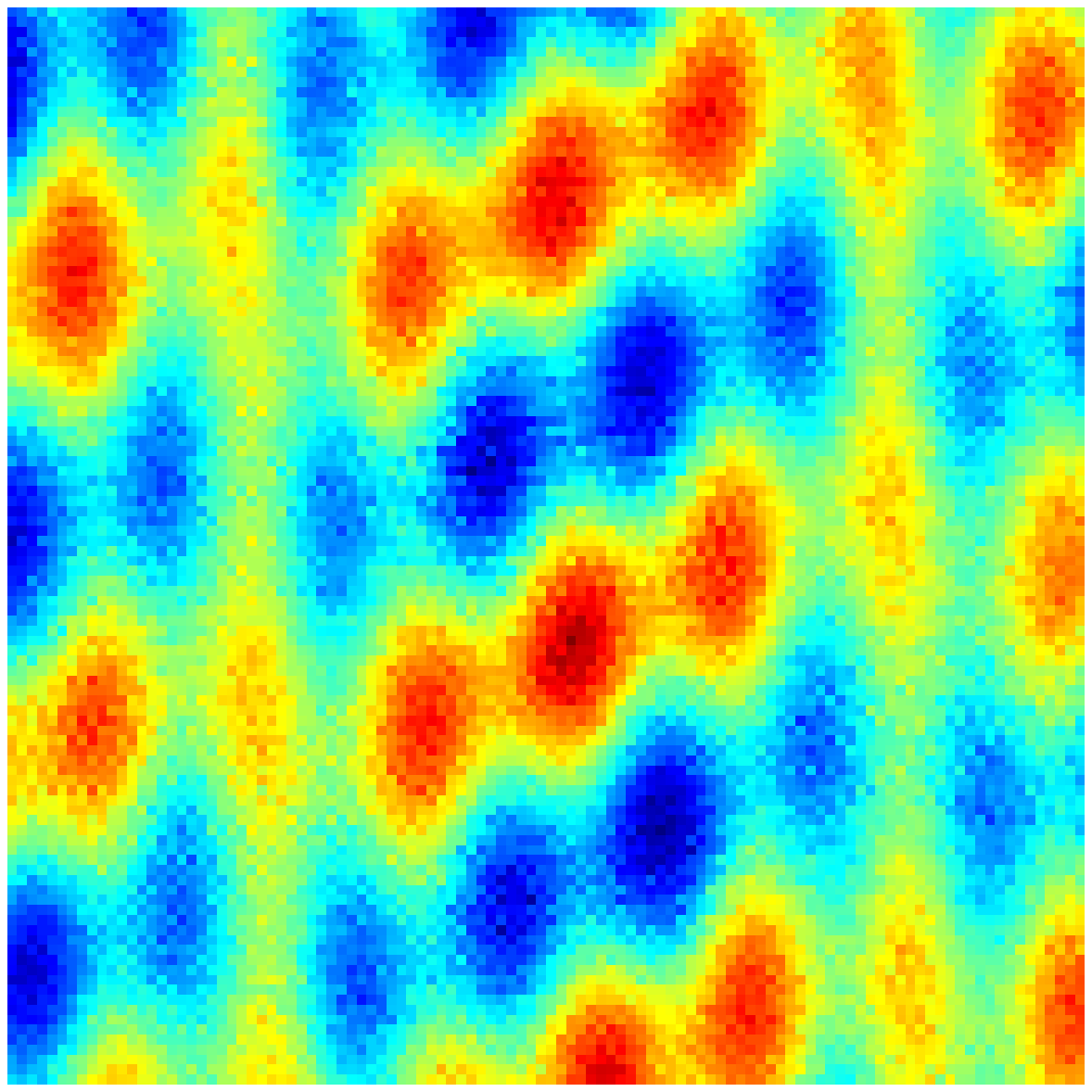}\\
\hbox{Standard recovery}&\hbox{Adaptive recovery}
\end{array} $$
\caption{Recovery for $\omega_{\max}=8.0$}
\label{pic:2}\end{figure}
\begin{figure}
$$\begin{array}{cc}
\epsfxsize=190pt\epsfysize=190pt\epsffile{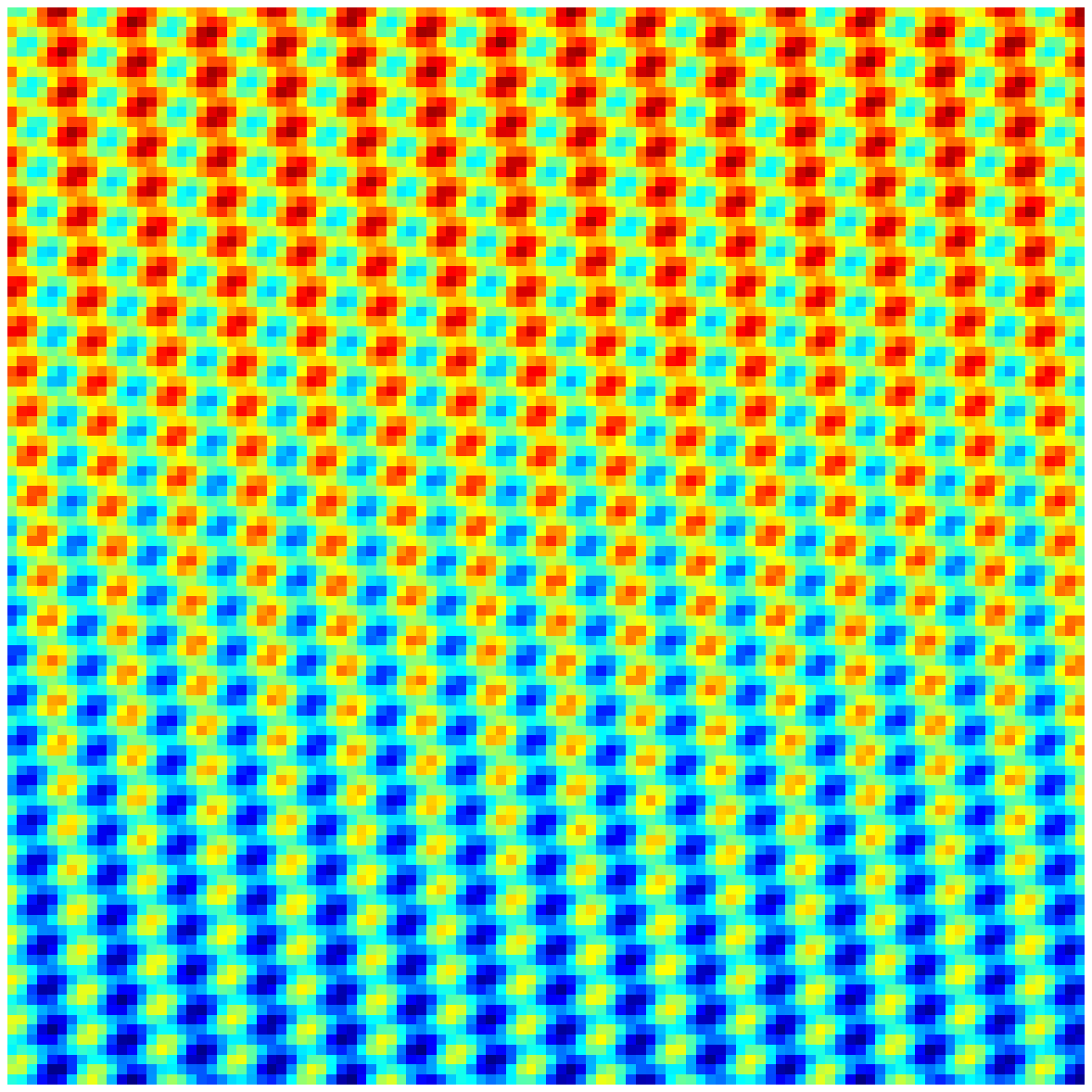}&
\epsfxsize=190pt\epsfysize=190pt\epsffile{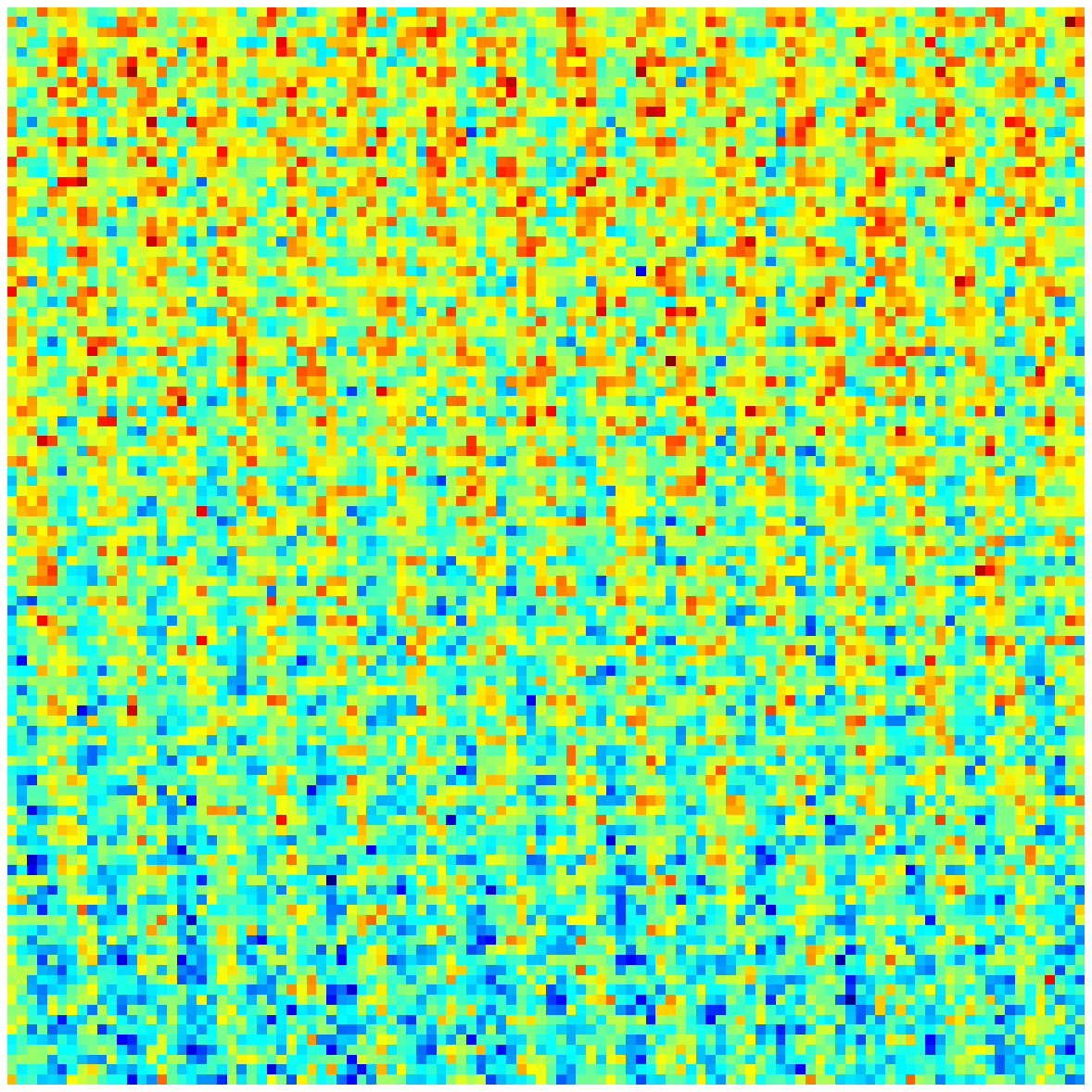}\\
\hbox{True Image}&\hbox{Observation} \\
\epsfxsize=190pt\epsfysize=190pt\epsffile{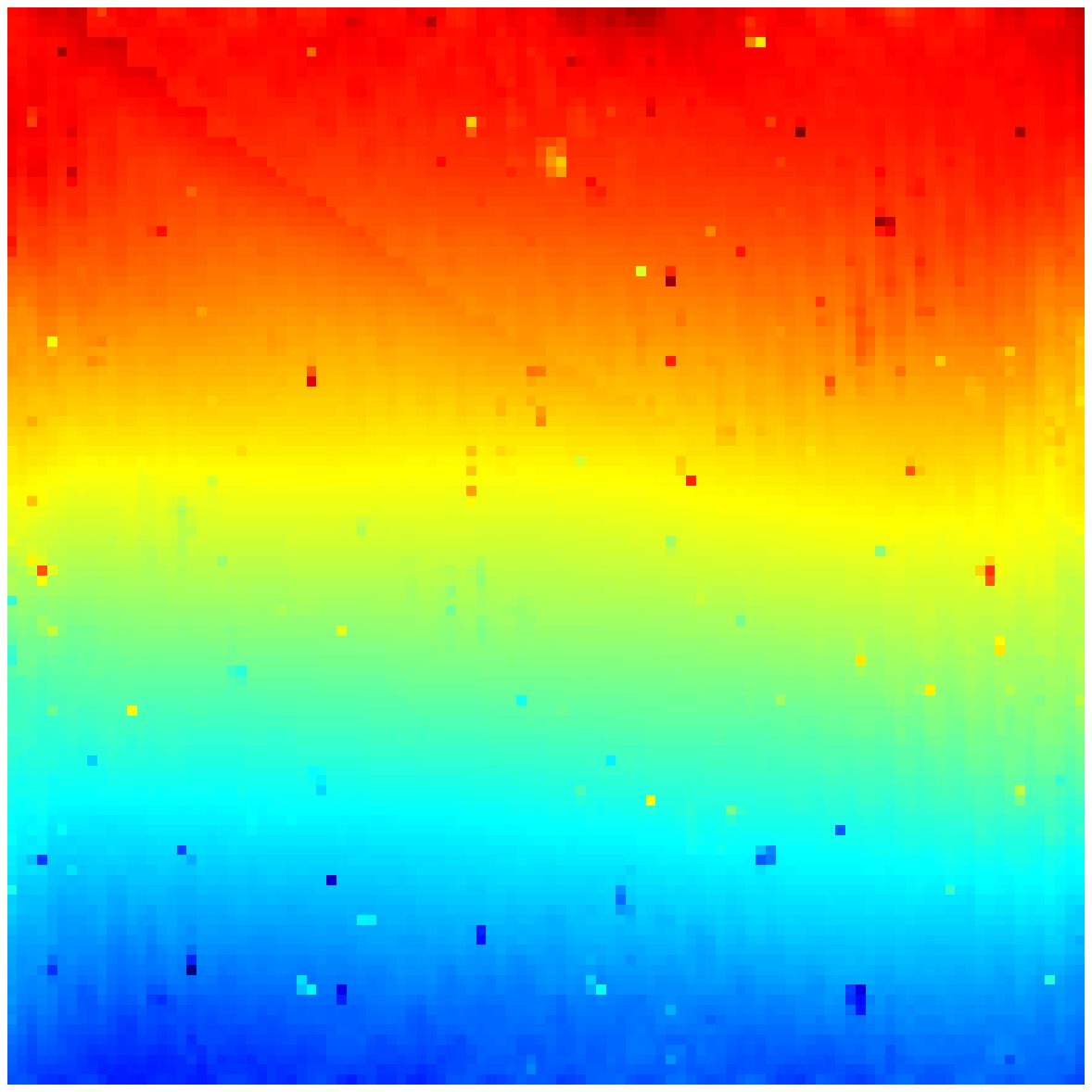}&
\epsfxsize=190pt\epsfysize=190pt\epsffile{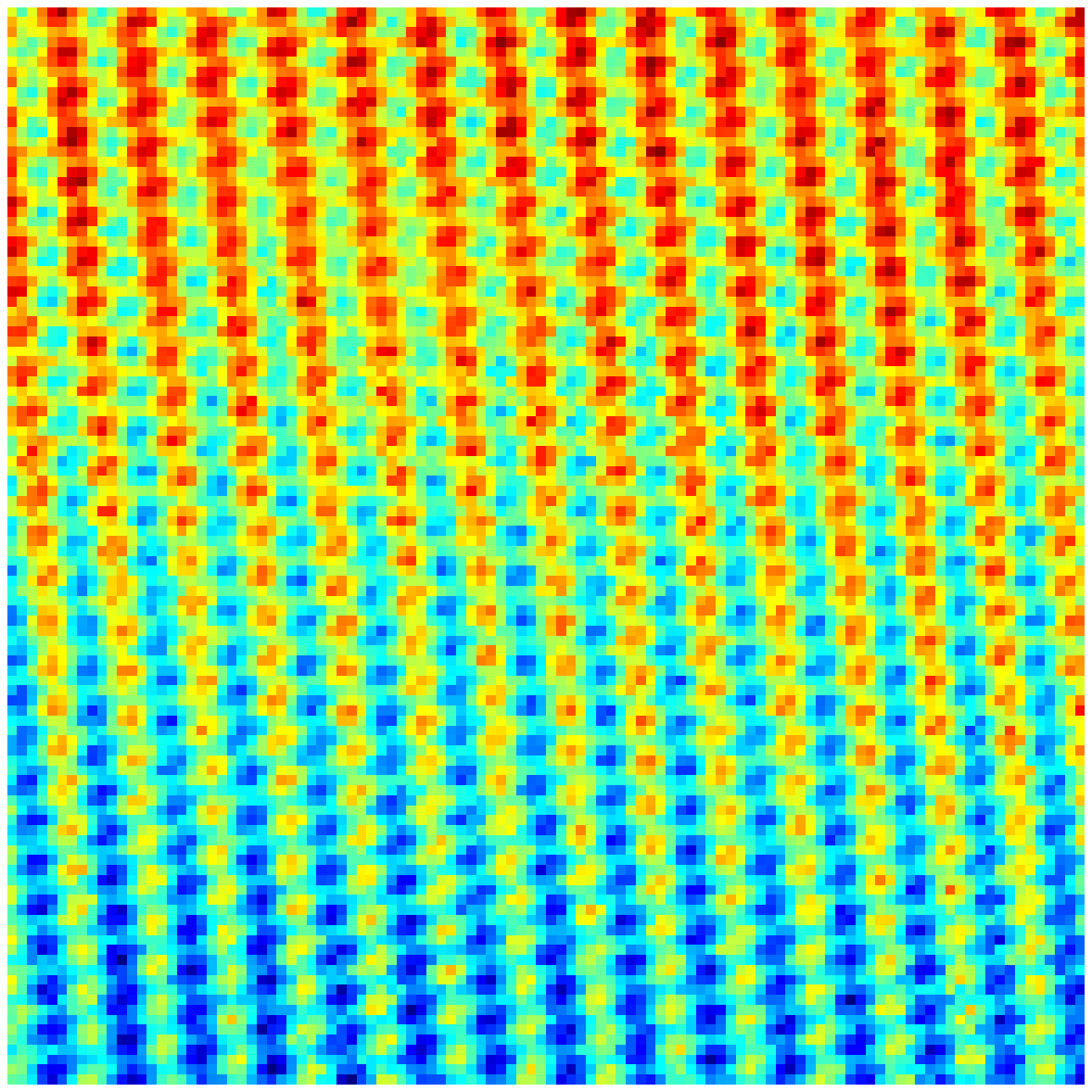}\\
\hbox{Standard recovery}&\hbox{Adaptive recovery}
\end{array} $$
\caption{Recovery for $\omega_{\max}=32.0$}
\label{pic:3}\end{figure}
\begin{figure}
$$\begin{array}{cc}
\epsfxsize=190pt\epsfysize=190pt\epsffile{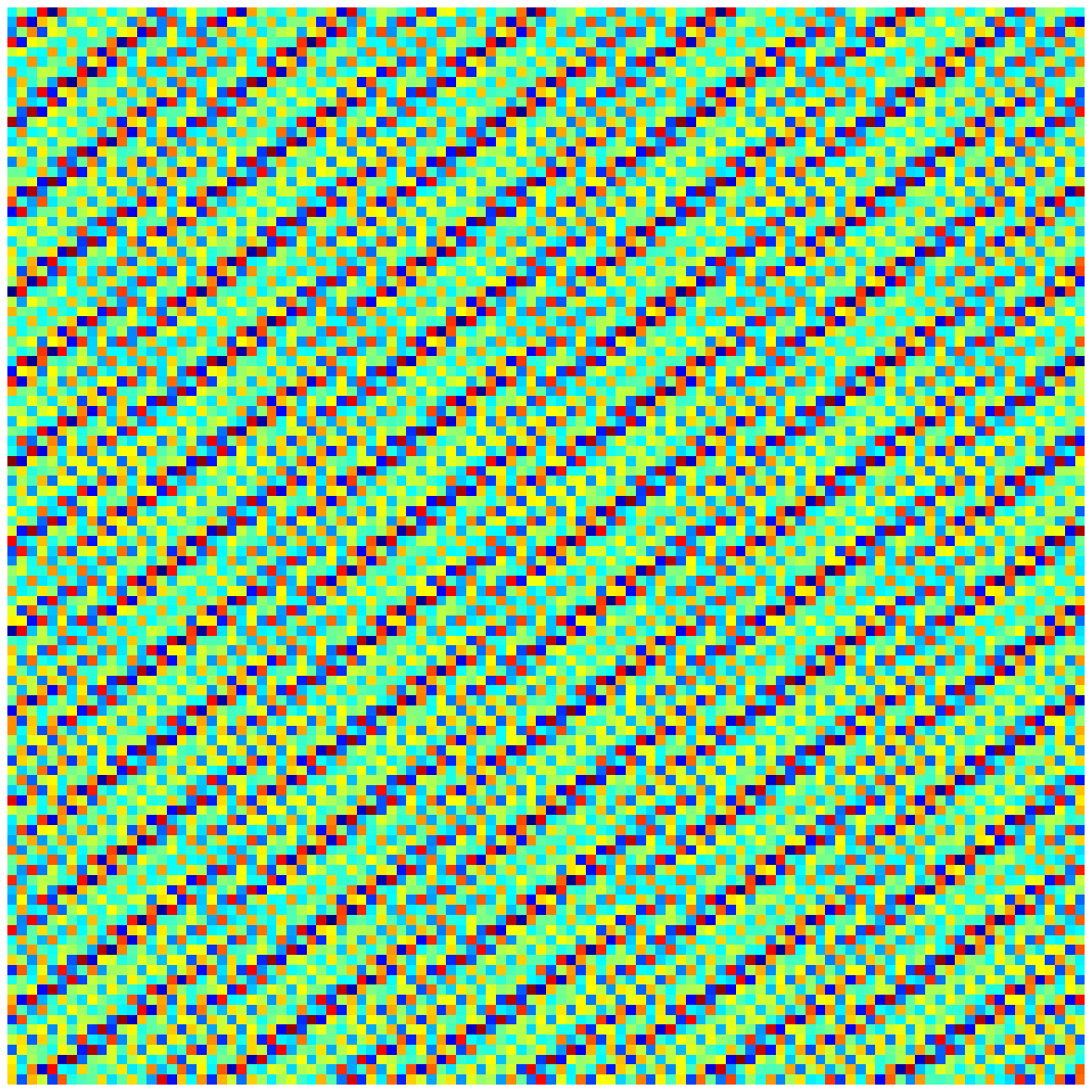}&
\epsfxsize=190pt\epsfysize=190pt\epsffile{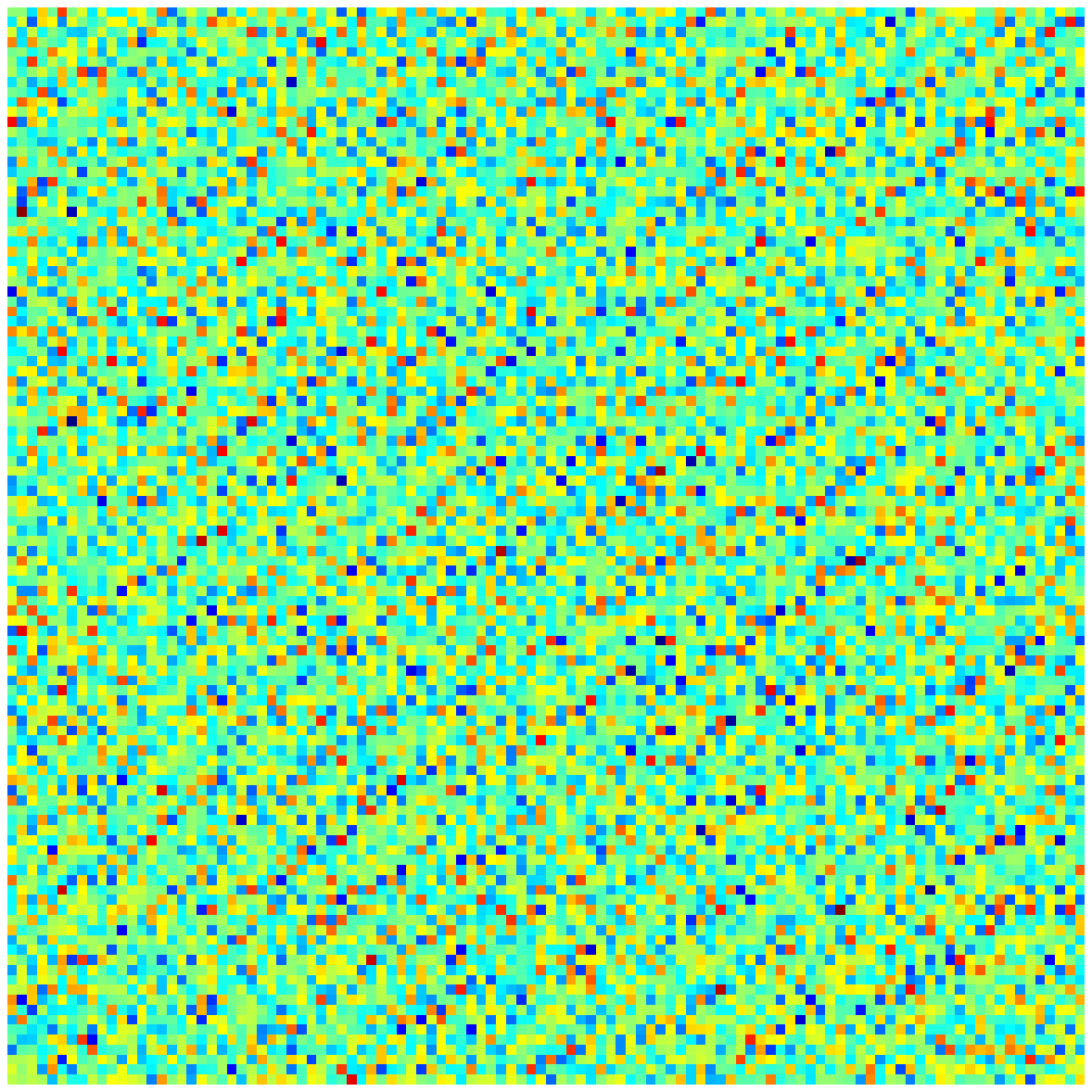}\\
\hbox{True Image}&\hbox{Observation} \\
\epsfxsize=190pt\epsfysize=190pt\epsffile{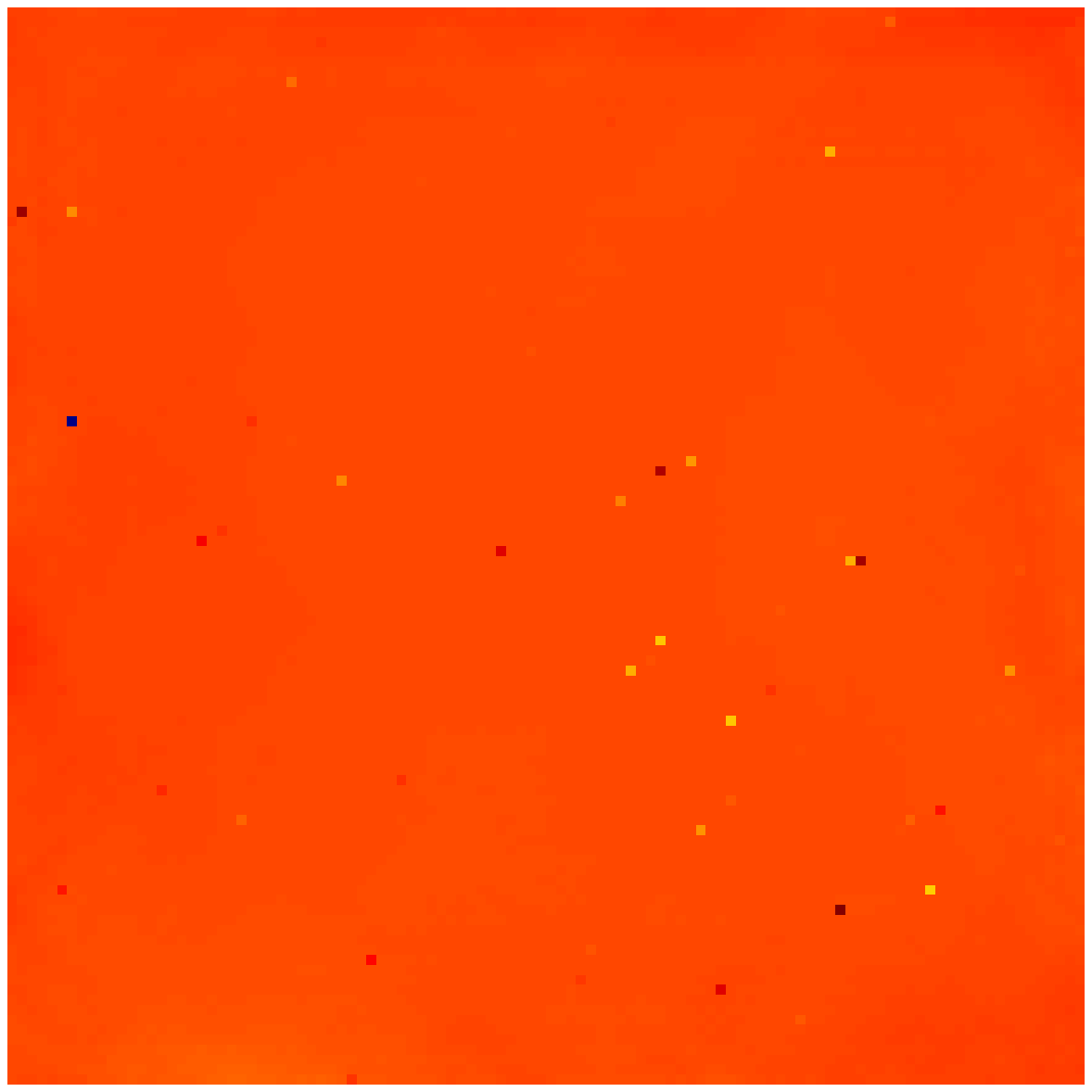}&
\epsfxsize=190pt\epsfysize=190pt\epsffile{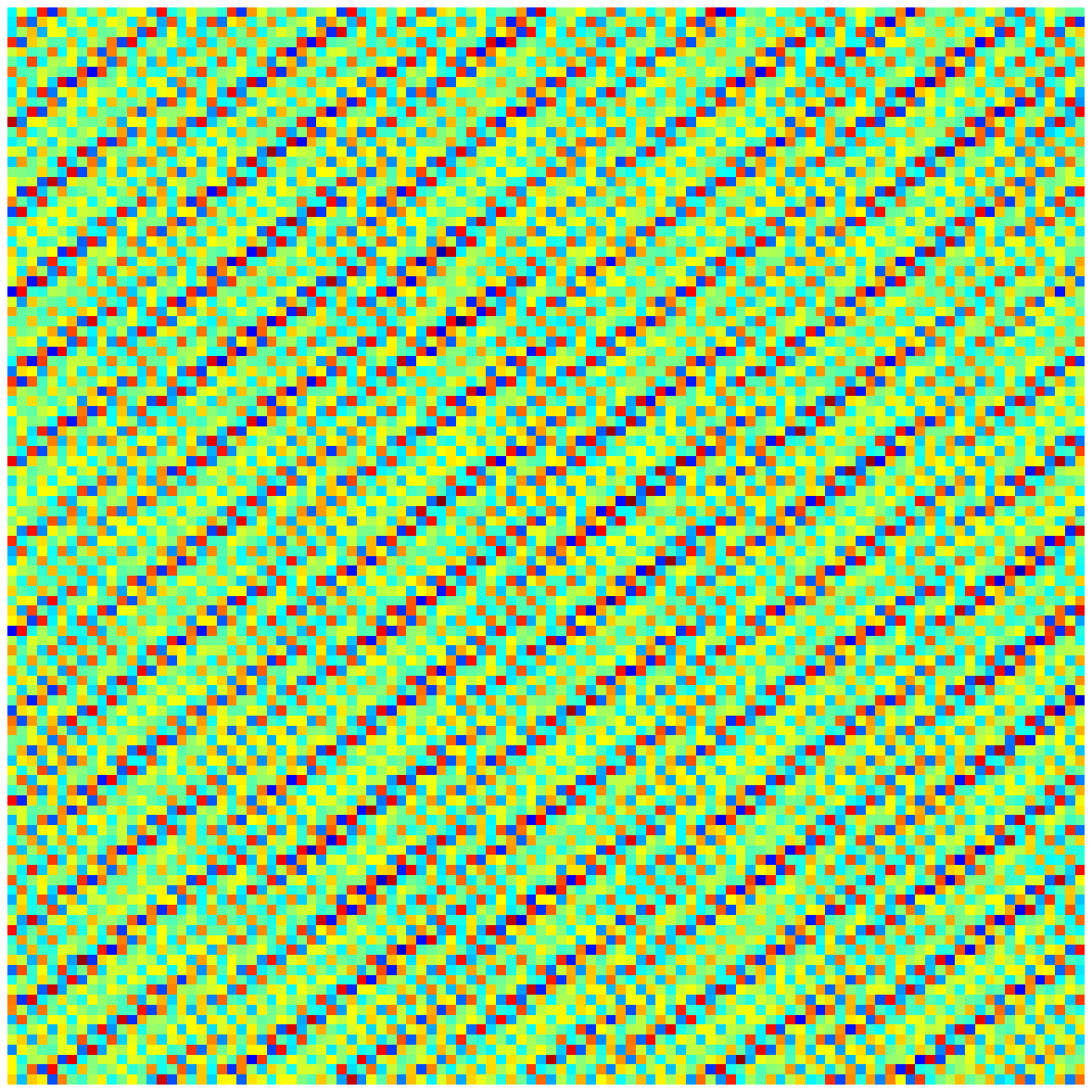}\\
\hbox{Standard recovery}&\hbox{Adaptive recovery}
\end{array} $$
\caption{Recovery for $\omega_{\max}=128.0$}\label{pic:4}\end{figure}

In the simulations we present here we compared the result of adaptive recovery with $T=10$ to that of a ``standard nonparametric recovery'', i.e. the recovery by the locally linear estimator with square window.
We have done $k=100$ independent runs for each of eight values of $\omega_{\max}$, \[\omega_{\max}=\{1.0,\,2.0,\,4.0,\,8.0\,16.0,\,32.0,\,64.0,\,128.0\}.
 \]In Table \ref{tab:1} we summarize the results for the mean integrated squared error (MISE) of the estimation,
\[
MISE=\sqrt{{1\over 100 m^2}\sum_{j=1}^{100}\sum_{\tau=(1,1)}^{(m,m)} (\widehat{s}^{(j)}_\tau-s^{(j)}_\tau)^2}.
\]
The observed phenomenon is rather expectable: for slowly oscillating signals the quality of the adaptive recovery is slightly worse than that of ``standard recovery'', which are tuned for estimation of regular signals. When we rise the frequency of the signal components, the adaptive recovery stably outperforms the standard recovery.
Finally, standard recovery is clearly unable to recover highly oscillating signals (cf Figures \ref{pic:1}-\ref{pic:4}).
\begin{table}[p]%
\caption{MISE of adaptive recovery}
\begin{tabular}{|c|c|c|}
\hline
& Standard& Adaptive\\[-1.5ex]
\raisebox{1.5ex}{$\omega_{\max}$} &  recovery &
     recovery \\
\hline\hline
$1.0$  &0.12   &0.1 \\\hline
$2.0$ &  0.20   &0.12   \\\hline
$4.0$&0.36  & 0.18 \\\hline
$8.0$& 0.54&0.27 \\\hline
$16.0$&   0.79 & 0.25 \\\hline
$32.0$&  0.75  & 0.29   \\\hline
$64.0$&     0.27   &  0.98 \\\hline
$128.0$&     0.24   &  1.00 \\
\hline
\end{tabular}
\label{tab:1}
\end{table}
\section*{Appendix}

We denote $C(\bZd)$ the linear space  of complex-valued fields over $\bZ^d$.
A field $r\in  C(\bZd)$ with finitely many nonzero entries
$r_{\tau}$ is called {\sl a filter}.
We use the commun notation  $\Delta_j$, $j=1,...,d$, for the ``basic shift
operators'' on $C(\bZd)$:
$$
(\Delta_jr)_{\tau_1,...,\tau_d}=r_{\tau_1,...,\tau_{j-1},
\tau_j-1,\tau_{j+1},...,\tau_d}.
$$
and denote $ r(\Delta)x$ the output of a filter  $r$, the
input to the filter being a field $x\in C(\bZd)$, so that $
(r(\Delta)x)_t=\sum\limits_\tau r_\tau x_{t-\tau}. $

\subsection{Proof of Lemma \ref{dapp:lemma1}.}
To save notation, let
$B=B_h(x)$ and $T=T_h(x)$. Let $p\in C(\bZd)$ be such that $p\in
{\cal S}^t_{3T}(0,\mu,T)$ and $|p_\tau-f_\tau|\leq
\Phi_\mu(f,{B}_h(x))$ for all $\tau\in \bZ(B_h(x))$. Since $p\in
{\cal S}^t_{3T}(0,\mu,T)$, there exists a filter $q\in C_T(\bZd)$
such that $|q|_2\leq \mu(2T+1)^{-d/2}$ and
$(q(\Delta)p)_\tau=p_\tau$ whenever $|\tau-t|\leq 3T$. Setting
$\delta_\tau=f_\tau-p_\tau$, we have for any $\tau$ , $|\tau-t|\leq 3T$,
\bse
\lefteqn{|f_\tau-(q(\Delta)f)_\tau|\leq|\delta_\tau|
+|p_\tau-(q(\Delta)p)_\tau|+|(q(\Delta)\delta)_\tau|}
\\
&\leq&\Phi_\mu(f,B_h(x))+|q|_1\max\{|\delta_\nu|:|\nu-\tau|\leq T\}
\leq\Phi_\mu(f,B_h(x))\\
&&+|q|_2(2T+1)^{d/2}\Phi_\mu(f,B_h(x))\max\{|\delta_\nu|:|\nu-\tau|\leq
T\}\\
&&{\hbox{[note that $|\tau-t|\leq 3T$ and
$|\nu-\tau|\leq T$ implies $|\nu-t|\leq 4T$]}}\\
&\leq&\Phi_\mu(f,B_h(x))(1+\mu)
={\Phi_\mu(f,B_h(x))(1+\mu)\over(2T+1)^{d/2}}(2T+1)^{-d/2}
\ese
as required in (\ref{dapp:em10}). \qed

 \subsection{Proof of Theorem \protect{\ref{3.themain}}}
 In the main body of the proof, we focus on the case
$p,q<\infty$; the case of infinite $p$ and/or $q$ will be
considered at the concluding step 5$^0$.
\\
Let us fix a family of well-filtered signals ${\cal F}={\cal
F}^{k,\rho,p}_d(B;R)$ with the parameters satisfying the premise
of Theorem \ref{3.themain} and a function $f$ from this class.
\\
Recall that by the definition of ${\cal F}$ there exists a
function $F\geq0$, $\|F\|_{p,B}\leq R$, such that for all
$x=m^{-1}t\in(\inter B)\cap\Gamma_n$ and all $h,\;B_h(x)\subset B$:
\begin{equation}
\label{3.peq1}
\Phi_\mu(f,{B}_h(x))\le P_1 h^{k-d/p} \Omega(f,{B}_h(x)),\,\,
\Omega(f,B')=\left(\displaystyle{\int\limits_{B'}} F^p(u)
du\right)^{1/p}.
\end{equation}
From now on, $P$ (perhaps with sub- or superscripts) are
quantities $\ge1$ depending on $\mu,d,\gamma,p$ only and
nonincreasing in $p>d$. \\~\par
{\bf 1$^0$.} We need the following auxiliary result:
\begin{lem}
\label{3.lem1} Assume that
\begin{equation}
\label{3.mbound1} n^{{{k-d/p}\over{d}}}\sqrt{\ln n} \ge
P_1(\mu+3)^{k-d/p+d/2}{R\over{\sigma\omega}}.
\end{equation}
Given a point $x\in \Gamma_n\cap B_\gamma$, let us choose the
largest $h=h(x)$ such that
\begin{equation}
\label{3.peq2}
\begin{array}{cl}
(a):&h\leq (1-\gamma)D(B),\\
(b):&P_1 h^{k-d/p}
\Omega(f,{B}_h(x))\leq S_n(h).
\end{array}
\end{equation}
Then $h(x)$ is well-defined and \begin{equation}\label{3.peq12}
h(x)\geq m^{-1}. \end{equation} Besides this, the error at $x$ of
the adaptive estimate $\widehat{f}_n$ as applied to $f$ can be
bounded as follows:
\begin{equation}
\label{3.peq3}
|\widehat{f}_n(x;y)-f(x)|\leq C_2 \left[S_n(h(x))\b1\{\epsilon\in\Xi_n\}+\sigma\Theta_{(n)}\b1\{\epsilon\not\in\Xi_n\}\right]
\end{equation}
\end{lem}
 {\bf Proof:} 
 The quantity $h(x)$ is
well-defined, since for small positive $h$ the left hand side in
(\ref{3.peq2}.$b$) is close to 0, while the right hand side one is
large. From (\ref{3.mbound1}) it
follows that $h=m^{-1}$ satisfies (\ref{3.peq2}.$a$), so that
${B}_{m^{-1}}(x)\subset B$. Moreover, (\ref{3.mbound1}.$b$) implies
that $$ P_1m^{-k+d/p}R\le S_n(m^{-1}); $$ the latter inequality, in
view of $\Omega(f,{B}_{m^{-1}}(x))\le R$, says that $h=m^{-1}$ satisfies
(\ref{3.peq2}.$b$) as well. Thus, $h(x)\ge m^{-1}$, as claimed in
(\ref{3.peq12}).
\par
Consider the window ${B}_{h(x)}(x)$. By
(\ref{3.peq2}.$a$) it is admissible for $x$, while from
(\ref{3.peq2}.$b$) combined with (\ref{3.peq1}) we get $
\Phi_\mu(f,{B}_{h(x)}(x))\le S_n(h). $ It follows that the ideal
window ${B}_*(x)$ of $x$ is not smaller
than ${B}_{h(x)}(x)$. 
\par Assume that $\epsilon\in\Xi_n$. Then, according to
(\ref{3.eq15}), we have
\begin{equation}
\label{3.peq4}
|\widehat{f}_n(x;y)-f(x)|\le5C_1\left[\Phi_\mu(f,{B}_{h_*(x)}(x))+S_n(h_*(x))\right].
\end{equation}
Now, by the definition of an ideal window,
 $\Phi_\mu(f,{B}_{h_*(x)}(x))\leq S_n(h_*(x))$,
and the right hand side in (\ref{3.peq4}) does not exceed $ 10C_1
S_n(h_*(x)) \le 10C_1S_n(h(x)) $ (recall that, as we have seen,
$h_*(x)\ge h(x)$), as required in
(\ref{3.peq3}). 
\par
Now let $\epsilon\not\in\Xi_n$. Note that
$\widehat{f}_n(x;y)$ is certain estimate $\widehat{f}^h(x;y)$
associated with a centered at $x$ and admissible for $x$ cube
${B}_h(x)$ which is normal and such that $h\geq m^{-1}$ (the
latter -- since the window $B_{m^{-1}}(x)$ always is normal, and
$B_{h}(x)$ is the largest normal window centered at $x$). Applying
(\ref{3.eq12}) with $h'=m^{-1}$ (so that
$\widehat{f}^{h'}_n(x;y)=f(x)+\sigma \epsilon_t$), we get $
|(f(x)+\sigma\epsilon_t)-\widehat{f}_n(x;y)|\leq 4C_1S_n(m^{-1}),
$ whence
$$
|f(x)-\widehat{f}_n(x;y)|\leq \sigma|\epsilon_t|+4C_1S_n(m^{-1})
\leq \sigma\Theta_{(n)}+4C_1\sigma\omega\sqrt{\ln n}\leq
C_2\Theta_{(n)}
$$
(recall that we are in the situation $\epsilon\not\in\Xi_n$,
whence $\omega\sqrt{\ln n}\leq \Theta_{(n)}$). We have arrived at
(\ref{3.peq3}). \qed \par
Now we are ready to complete the proof. Assume that
(\ref{3.mbound1}) takes place, and let us fix $q$, ${{2k+d}\over
d}p\le q<\infty$.
\medskip
\\
{\bf 2$^0$}. Let us denote $\widehat{\sigma}_n=\sigma\sqrt{{{\ln n}\over n}}$. Note that for every $x\in \Gamma_n\cap B_\gamma$ either
$$
h(x)=(1-\gamma)D(B),
$$
or
$$
h(x)=\left({{\widehat{\sigma}_n}\over{P_1\Omega(f,B_{h(x)}(x))}}\right)^{{2p\over{2kp+(p-2)d}}},
$$
what means that
\begin{equation}\label{3.nnew1}
P_1h^{k-d/p}(x)\Omega(f,B_{h(x)}(x))=
S_n(h(x)).
\end{equation}
Let $U,V$ be the sets of those $x\in B_\gamma^n\equiv \Gamma_n\cap
B_\gamma$ for which the first or, respectively, the second of this
possibilities takes place. If $V$ is nonempty, let us partition it
as follows.
\\
 1) We can choose $x_1\in V$ ($V$ is finite!) such that
$ h(x)\ge h(x_1)\quad \forall x\in V. $ After $x_1$ is chosen, we
set $ V_1=\{x\in V\mid\, B_{h(x)}(x)\cap
B_{h(x_1)}(x_1)\ne\emptyset\}. $
\\
2) If the set $V\backslash V_1$ is nonempty, we apply the
construction from 1) to this set, thus getting $x_2\in V\backslash
V_1$ such that $ h(x) \ge h(x_2) \quad \forall x\in V\backslash
V_1, $ and set $ V_2=\{x\in V\backslash V_1\mid\, B_{h(x)}(x)\cap
B_{h(x_2)}(x_2)\ne\emptyset\}. $ If the set $V\backslash(V_1\cup
V_2)$ still is nonempty, we apply the same construction to this
set, thus getting $x_3$ and $V_3$, and so on.\\ The outlined
process clearly terminates after certain step (since $V$ is
finite).   On termination, we get a collection of $M$ points
$x_1,...,x_M\in V$ and a partition $ V=V_1\cup V_2\cup...\cup V_M
$ with the following properties:
\begin{enumerate}
\item  The cubes $B_{h(x_1)}(x_1),...,B_{h(x_M)}(x_M)$ are mutually
disjoint;
\item For every $\ell\le M$ and every $x\in V_\ell$ we have $
h(x)\ge h(x_\ell)\hbox{\ and\ } B_{h(x)}(x)\cap
B_{h(x_\ell)}(x_\ell) \ne\emptyset. $
\\
We claim that also
\item For every $\ell\le M$ and every $x\in V_\ell$ one has
\begin{equation}
\label{3.nnew2} h(x)\ge
\max\left[h(x_\ell);\|x-x_\ell\|_\infty\right].
\end{equation}
\end{enumerate}
Indeed, $h(x)\ge h(x_\ell)$ by (ii), so that it suffices to verify
(\ref{3.nnew2}) in the case when $\|x-x_\ell\|_\infty\ge
h(x_\ell)$. Since $B_{h(x)}(x)$ intersects
$B_{h(x_\ell)}(x_\ell)$, we have
$$\|x-x_\ell\|_\infty\le
{1\over2}(h(x)+h(x_\ell)).$$
Whence
$$ h(x) \ge
2\|x-x_\ell\|_\infty-h(x_\ell)\ge \|x-x_\ell\|_\infty, $$
which is what we need. \\
{\bf 3$^0$.} Let us set $B_\gamma^n=\Gamma_n\cap B_\gamma$.
Assume that $\epsilon\in\Xi_n$. When substituting  $h(x)=(1-\gamma)[D(B)]$ for $x\in U$, we have by (\ref{3.peq3}):
\bse
\lefteqn{\left|\widehat{f}_n(\cdot;y)-f(\cdot)\right|_{q,B_\gamma}^{q}\leq C_2^{q}m^{-{d\over
q}}{\displaystyle}{\sum\limits_{x\in B_\gamma^n}} S_n^{q}(h(x))}\\
&=&C_2^{q}m^{-{d\over q}}{\displaystyle}{\sum\limits_{x\in
U}}S_n^{q}(h(x))
+C_2^qm^{-{d\over
q}}{\displaystyle}{\sum\limits_{\ell=1}^M\sum\limits_{x\in
V_\ell}} S_n^{q}(h(x))\\
&=&C_2^{q}m^{-{d\over
q}}{\displaystyle}{\sum\limits_{x\in U}
\left[{{\widehat{\sigma}_n}\over{((1-\gamma)D(B))^{d/2}}}\right]^q}
+C_2^qm^{-{d\over
q}}{\displaystyle}{\sum\limits_{\ell=1}^M\sum\limits_{x\in V_\ell}
S_n^{q}(h(x))} \\
\mbox{[by (\ref{3.nnew2})]}\;&\leq& C_3^q\widehat{\sigma}_n^q
m^{-{d\over q}}{\displaystyle}{\sum\limits_{\ell=1}^M
\sum\limits_{x\in V_\ell}
\left(\max\left[h(x_\ell),\|x-x_\ell\|_\infty\right]\right)^{-{dq\over2}}}
+C_3^q\widehat{\sigma}_n^q[D(B)]^{{d(2-q)\over 2}}\\
&\leq&
C_4^q\widehat{\sigma}_n^q
{\displaystyle}{\sum\limits_{\ell=1}^M \int
\left(\max\left[h(x_\ell),\|x-x_\ell\|_\infty\right]\right)^{-{dq\over2}}}dx
+C_3^q\widehat{\sigma}_n^q[D(B)]^{{d(2-q)\over 2}} \\
&\leq&
C_5^q\widehat{\sigma}_n^q {\displaystyle}{\sum\limits_{\ell=1}^M
\int\limits_{0}^\infty r^{d-1}
\left(\max\left[h(x_\ell),r\right]\right)^{-{dq\over2}}dr}+C_3^q\widehat{\sigma}_n^qD[D(B)]^{{d(2-q)\over
2}},
\ese
due to $h(x)\geq m^{-1}$, see (\ref{3.peq12}). Further,
note that $${dq\over2}-d+1\ge {{2k+d}\over 2}p-d+1\ge d^2/2+1$$
in view of $q\ge {{2k+d}\over d}p$, $k\ge1$ and $p>d$, and
\bse
\lefteqn{\left|\widehat{f}_n(\cdot;y)-f(\cdot)\right|_{q,B_\gamma}^{q}
\leq C_6^q\widehat{\sigma}_n^q
{\displaystyle}{\sum\limits_{\ell=1}^M
\left[h(x_\ell)\right]^{{d(2-q)\over
2}}}+C_3^q\widehat{\sigma}_n^q[D(B)]^{{d(2-q)\over 2}}}\\
\mbox{[by (\ref{3.nnew1}]]}\;
&=& C_6^q\widehat{\sigma}_n^q
{\displaystyle}{\sum\limits_{\ell=1}^M}
\left[{{\widehat{\sigma}_n}\over{P_1
\Omega(f,B_{h(x_\ell)}(x_\ell))}}\right]^{{d(2-q)\over2k-2d/p+d}}
+C_3^q\widehat{\sigma}_n^q[D(B)]^{{d(2-q)\over 2}}\\
&=&C_3^q\widehat{\sigma}_n^q[D(B)]^{{d(2-q)\over 2}}
+C_6^q\widehat{\sigma}_n^{2\beta(p,k,d,q)q}
{\displaystyle}{\sum\limits_{\ell=1}^M}
\left[P_1\Omega(f,B_{h(x_\ell)}(x_\ell))\right]^{{d(q-2)\over2k-2d/p+d}}
\ese
 by definition of
$\beta(p,k,d,q)$.
\\
Now note that ${{d(q-2)}\over{2k-2d/p+d}}\ge p$ in view of $q\ge
{{2k+d}\over d}p$, so that
\bse
\sum\limits_{\ell=1}^M
\left[P_1\Omega(f,B_{h(x_\ell)}(x_\ell))\right]^{{{d(q-2)}\over{2k-2d/p+d}}}
&\leq&\left[\sum\limits_{\ell=1}^M
\left(P_1\Omega(f,B_{h(x_\ell)}(x_\ell))\right)^p\right]^{{{dq-2d}\over{p(2k-2d/p+d)}}}
\\
&\leq&\left[P_1^pR^p\right]^{{{d(q-2)}\over{p(2k-2d/p+d)}}}
\ese
(see (\ref{3.peq1}) and take into account that the cubes
$B_{h(x_\ell)}(x_\ell)$, $\ell=1,...,M$, are mutually disjoint by
(i)). We conclude that for $\epsilon\in\Xi_n$
\begin{eqnarray}
\label{3.nnew9}
\left|\widehat{f}_n(\cdot;y_f(\epsilon))-f(\cdot)\right|_{q,B_\gamma}&\le&
C_7\widehat{\sigma}_n[D(B)]^{{d(2/q-1)\over2}}+
P_2\widehat{\sigma}_n^{2\beta(p,k,d,q)}
R^{{{d(1-2/q)}\over{2k-2d/p+d}}}\nonumber \\
&=&C_7\widehat{\sigma}_n[D(B)]^{{d(2/q-1)\over2}} +
P_2R\left({{\widehat{\sigma}_n}\over R}\right)^{2\beta(p,k,d,q)}.
\end{eqnarray}
{\bf 4$^0$.} Now assume that $\epsilon\not\in\Xi_n$. In this
case, by (\ref{3.peq3}), $$ |\widehat{f}_n(x;y)-f(x)|\le
C_2\sigma\Theta_{(n)}\quad \forall x\in B_\gamma^n. $$ Hence,
taking into account that $mD(B)\geq1$,
\begin{equation}
\label{3.peq23}
\left|\widehat{f}_n(\cdot;y)-f(\cdot)\right|_{q,B_\gamma}
\le C_2\sigma\Theta_{(n)} [D(B)]^{{d\over q}}.
\end{equation}
{\bf 5$^0$.}
When combining (\ref{3.nnew9}) and (\ref{3.peq23}), we
get
$$
 \left(E
\left\{\|\widehat{f}_n(\cdot;y)-f(\cdot)\|_{q,B_\gamma}^2\right\}\right)^{1/2}\\
\leq C_8\max\left[\widehat{\sigma}_n[D(B)]^{{{d(2/q-1)}\over2}};\,
P_4
 R\left({{\widehat{\sigma}_n}\over{R}}\right)^{2\beta(p,k,d,q)};\;J\right],
 $$
 where
 \bse
 J^2&=&E \left\{\b1\{\epsilon\not\in \Xi_n\}
C_2\sigma^2\Theta_{(n)}^2\right\} \le  C^2_2\sigma^2
P^{{1/2}}\{\epsilon\not\in \Xi_n\}\left(E \left\{\Theta_{(n)}^4\right\}\right)^{{1/2}}\\
&\le&  C_9\sigma^2 n^{-2(\mu+1)}{\ln n}
\ese
(we have used (\ref{3.beq18}) and (\ref{3.beq18a})). Thus, when (\ref{3.mbound1}) holds, for all $d<p<\infty$ and all $q$,
${{2k+d}\over d} p\le q<\infty$ we have
\begin{eqnarray}
\label{3.peqfin}
\lefteqn{\left(E
\left\{\|\widehat{f}_n(\cdot;y)-f(\cdot)\|_{q,B_\gamma}^2\right\}\right)^{1/2}}\nonumber\\
&\le& C_8\max\left[ \widehat{\sigma}_n[D(B)]^{{{d(2/q-1)}\over2}},
 P_4
 R\left({{\widehat{\sigma}_n}\over{R}}\right)^{2\beta(p,k,d,q)},
{C^{1/2}_9\sigma \sqrt{\ln n}\over n^{(\mu+1)}}\right].
\end{eqnarray}
Now it is easily seen that if $P\ge1$ is a properly chosen
function of $\mu,d,\gamma,p$ nonincreasing in $p>d$ and
(\ref{3.large}) takes place then
\begin{enumerate}
\item assumption
(\ref{3.mbound1}) holds,
\item
the right hand side in
(\ref{3.peqfin}) does not exceed the quantity
$$
PR\left({{\widehat{\sigma}_n}\over{R}}\right)^{2\beta(p,k,d,q)}=
PR\left({{\widehat{\sigma}_n}\over{R}}\right)^{2\beta(p,k,d,q)}
[D(B)]^{d\lambda(p,k,d,q)}
$$
(recall that $q\ge {{2k+d}\over d}p$, so that $\lambda(p,k,d,q)=0$).
\end{enumerate}
We conclude the bound (\ref{3.main}) for the case of
$d<p<\infty$, $\infty> q\ge {{2k+d}\over d}p$. When passing to the limit as
$q\to\infty$, we get the desired bound for $q=\infty$ as well.
\par
Now let $d<p<\infty$ and $1\le q\le q_*\equiv
{{2k+d}\over d}p.$ By the H\"older inequality and in view of
$mD(B)\geq1$ we have $$ \left|g\right|_{q,B_\gamma} \le
C_{10}\left|g\right|_{q_*,B_\gamma}|B_\gamma|^{{1\over q}-{1\over{q_*}}}, $$
and thus $$ \widehat{{\cal R}}_q\left(\widehat{f}_n;{\cal F}\right)
\le C_{10}\widehat{{\cal R}}_{q_*}\left(\widehat{f}_n;{\cal
F}\right)[D(B)]^{d\left({1\over q}-{1\over q_*}\right)}. $$
Combining
this observation with the (already proved) bound (\ref{3.main})
associated with $q=q_*$, we see that (\ref{3.main}) is valid for
all $q\in[1,\infty]$, provided that $d<p<\infty$. Passing in the
resulting bound to limit as $p\to\infty$, we conclude the validity of
(\ref{3.main}) for all $p\in(d,\infty]$,
$q\in[1,\infty]$. \qed

\newpage
\newpage

\newpage
\newpage

\end{document}